\documentclass[12pt]{article}

\usepackage{graphicx}
\usepackage{amsmath,amssymb}
\usepackage{epstopdf}
\usepackage{float}
\usepackage[small,nohug,heads=littlevee]{diagrams} \diagramstyle[labelstyle=\scriptstyle] 
\numberwithin{equation}{section}

\begin{document}

\begin{flushright}{YITP-SB-07-39}\end{flushright}

\begin{center}
{\LARGE Elliptic constructions of hyperk\"ahler metrics I: \\[5pt] The Atiyah-Hitchin manifold}
\vskip40pt
{\bf Radu A. Iona\c{s}}
\end{center}
\vskip10pt
\centerline{\it C.N.Yang Institute for Theoretical Physics, Stony Brook University}
\centerline{\it Stony Brook, NY 11794-3840, USA }
\centerline{\tt ionas@max2.physics.sunysb.edu}
\vskip60pt

\begin{abstract}
\noindent This is the first in a series of papers in which we develop a twistor-based method of constructing hyperk\"ahler metrics from holomorphic functions and elliptic curves. As an application, we revisit the Atiyah-Hitchin manifold and derive in an explicit holomorphic coordinate basis closed-form formulas for, among other things, the metric, the holomorphic symplectic form and all three K\"ahler potentials.
\end{abstract}

\newpage

\tableofcontents

\vspace*{40pt}

\setcounter{section}{-1}

\section{Introduction}

The generalized Legendre transform (GLT) approach to constructing hyperk\"ahler metrics emerged originally in the context of ${\cal N}=2$ supersymmetric sigma-models as the by-product of a superspace generalization of the Hodge duality between 0 and 2-form gauge fields in 4 dimensions \cite{Lindstrom:1983rt}. The close connection between supersymmetric sigma-models and complex geometry became clear in \cite{Zumino:1979et} and in particular, the target space geometry of 4-dimensional sigma-models with ${\cal N}=2$ supersymmetry was shown in \cite{AlvarezGaume:1981hm} to be hyperk\"ahler.

The natural mathematical setting of Lindstr\"om and Ro\v{c}ek's construction was recognized to be the theory of twistor spaces of hyperk\"ahler manifolds \cite{Hitchin:1986ea}. This is a classic  generalization of Penrose's non-linear graviton construction \cite{MR0439004} and forms an integral part of Salamon's theory of twistor spaces of quaternionic-K\"ahler manifolds \cite{MR664330}. Every hyperk\"ahler manifold possesses a 2-sphere's worth of holomorphic symplectic structures. The twistor space $Z$ of a hyperk\"ahler manifold ${\cal M}$ can be viewed as a holomorphic fibration $Z \rightarrow \mathbb{CP}^1$, with the fiber over a generic point $\zeta \in \mathbb{CP}^1 \simeq S^2$ being a copy of ${\cal M}$ endowed with the holomorphic symplectic structure corresponding to $\zeta$. In the twistor space picture the GLT equations arise through a patching construction of a twisted holomorphic symplectic bundle over $Z$ by means of twisted canonical transformations of type II. The metric information of the hyperk\"ahler manifold is encoded in a single holomorphic function of one or several sections of ${\cal O}(2j)$ bundles over $Z$ that satisfy a reality condition with respect to the real structure induced on $Z$ by antipodal conjugation on the sphere of complex structures. The Penrose transform of this function is a real function $F$ of the parameters of the ${\cal O}(2j)$ sections which satisfies a set of second order differential equations. A generalized Legendre-Fenchel transform of $F$ yields both a K\"ahler potential of the hyperk\"ahler manifold and a corresponding set of holomorphic coordinates. In this coordinate basis the components of the metric are given by the second derivatives of $F$ and take the form of ratios of determinants of Hankel matrices. 

An alternative route to the GLT was proposed in \cite{MR1848654}. It relies on the notion of twistor group and explains the GLT equations through a generalized hyperk\"ahler quotient construction.

Many if not all the metrics constructed by means of the GLT based on ${\cal O}(2j)$ sections can be interpreted as metrics on moduli spaces of monopoles. Two main approaches to constructing monopoles are used in the literature: Hitchin's twistor approach \cite{MR649818,MR709461} and Nahm's equations \cite{MR807414}. In the twistor theory, a $k$-monopole is equivalent to a degree $k$ spectral curve on the tangent bundle of $\mathbb{CP}^1$ satisfying a set of conditions. Among other things, a certain line bundle over the curve is required to be trivial. This condition was shown by Ercolani and Sinha to be equivalent to an algebraic constraint \cite{Ercolani:1989tp}. The moduli spaces of $k$-monopoles are coordinatized in this frame by means of degree $k$ based rational maps between Riemann spheres \cite{MR769355}. On the other hand, Nahm's approach is an adaptation of the ADHM construction of instantons, replacing matrices by differential operators and the quadratic constraints on the matrices by non-linear ODE's for three $k \times k$ matrices. The solutions of these differential equations with certain boundary and regularity conditions correspond to monopole solutions of the Bogomol'nyi equations. In other words, the moduli spaces of solutions are the same.

${\cal O}(2)$-based generalized Legendre transform constructions of hyperk\"ahler metrics have been extensively discussed in the literature. These constructions generalize the Gibbons-Hawking class of 4-dimensional hyperk\"ahler metrics with a locally toric action \cite{Gibbons:1979zt}, giving $4n$-dimensional hyperk\"ahler metrics with $n$ commuting isometries \cite{Hitchin:1986ea,MR953820}.  Prominent examples include the asymptotically locally Euclidean (ALE) and asymptotically locally flat (ALF) metrics of type $A_n$ (multi-Eguchi-Hanson and multi-Taub-NUT) \cite{Gibbons:1979zt} and the 8-dimensional Swann bundle used by Calderbank and Pedersen to classify  selfdual Einstein metrics with two commuting isometries \cite{Calderbank:2001uz} and by Anguelova, Ro\v{c}ek and Vandoren to describe the geometry of the classical moduli space of the universal hypermultiplet in string theory compactifications \cite{Anguelova:2004sj}.

Comparatively, the study of ${\cal O}(2j)$-based constructions with $j\geq2$ is less developed, due to inherent difficulties. Explicit forms of the metrics are known only in the rare cases when extraneous, usually symmetry-based arguments could be used to construct them. On the other hand, we have a much better understanding of the GLT holomorphic potentials. The 2-monopole moduli space metric of Atiyah and Hitchin approaches asymptotically, up to exponentially suppressed corrections, a Taub-NUT metric with negative mass parameter \cite{Atiyah:1988jp}. The latter metric has a GLT description in terms of a real ${\cal O}(2)$ section \cite{Karlhede:1984vr}. This fact was exploited by Ivanov and Ro\v{c}ek to conjecture the form of the GLT holomorphic potential describing the Atiyah-Hitchin metric \cite{Ivanov:1995cy}. Specifically, this was obtained from the Taub-NUT potential by essentially replacing the ${\cal O}(2)$ section with the square root of an ${\cal O}(4)$ one. Later on, Houghton gave a rigorous justification of this procedure and generalized it in fact to monopoles with charge $k$ within the frame of Hitchin's twistor approach to monopoles \cite{Houghton:1999hr}. In a similar manner, the ${\cal O}(2)$-based GLT description of ALE and ALF spaces of type $A_n$ was used by Chalmers, Ro\v{c}ek and Wiles to put forward and test a conjectural holomorphic potential corresponding to  the type $D_n$ spaces \cite{Chalmers:1998pu}. Rigorous justifications along the lines of Nahm's and Hitchin's approaches to monopoles were provided subsequently by Cherkis, Hitchin and Kapustin \cite{Cherkis:1998xca,Cherkis:2003wk}. 

In this and the following papers of this series \cite{MM2,MM3}, we investigate in detail the ${\cal O}(4)$-based GLT constructions  and develop in the process a generic approach that allows one to perform calculations explicit to a high degree.

In sections \ref{SEC:twist-sp-hk} through \ref{SEC:hk-metrics-GLT} we review the theory of twistor spaces of hyperk\"ahler manifolds and the GLT approach to constructing hyperk\"ahler metrics, following mostly \cite{Hitchin:1986ea,Lindstrom:1987ks}. In section \ref{SEC:O2-GLT} we review the ${\cal O}(2)$-based constructions and, in particular, the construction of the Taub-NUT metric in the frame of the GLT, which was done for the first time in \cite{Karlhede:1984vr}. In section \ref{SEC:O4-multiplets} we begin a generic analysis of the 2-monopole spectral curve. In section \ref{SEC:AH-metric} we complete the program started in \cite{Ivanov:1995cy} and apply the machinery developed in the previous section as well as in the Appendix (which is, for the most part, and advertisement for the use of Jacobi's elliptic integrals of the third kind in Weierstrass form) to construct explicitly the geometry of the Atiyah-Hitchin metric.

\section{Twistor spaces of hyperk\"ahler manifolds} \label{SEC:twist-sp-hk}

We begin with a review of the theory of twistor spaces of hyperk\"ahler manifolds \cite{Hitchin:1986ea}. Let ${\cal M}$ be a hyperk\"ahler manifold, with the standard triplet of complex structures $J_1$, $J_2$, $J_3$ forming a quaternionic algebra
\begin{equation}
J_1^2 = J_2^2 = J_3^2 = - I\, , \qquad J_1J_2 = J_3\, , \qquad \mbox{a.s.o.} \label{quater-neon}
\end{equation}
where $I$ is the identity endomorphism on the tangent bundle of ${\cal M}$. These generate in fact a 2-sphere's worth of complex structures on ${\cal M}$: for any unit vector $(x^1,x^2,x^3) \in \mathbb{R}^3$, the linear combination $x^1J_1+x^2J_2+x^3J_3$ forms an integrable complex structure on ${\cal M}$, compatible with the metric and Levi-Civita connection. In particular, the quaternionic algebra properties of $J_1$, $J_2$, $J_3$ guarantee that $(x^1J_1+x^2J_2+x^3J_3)^2 = - I$. Identifying the 2-sphere with the complex projective space $\mathbb{CP}^1$, one can  describe it alternatively in terms of a set of two homogeneous coordinates, which we denote here by $\pi^A$, with  $A\in\{1,2\}$. The projective and the extrinsic descriptions of the Riemann sphere are related through the stereographic projection
\begin{equation}
x^i =  \frac{\bar{\pi}_A(\sigma^i)^A{}_B\pi^B}{\bar{\pi}_C\pi^C} \label{stereo}
\end{equation}
where $i \in \{1,2,3\}$ is a Cartesian index, $\sigma^i$ are the three standard $2\times 2$ Pauli matrices and $\bar{\pi}_A$ denotes the complex conjugate of $\pi^A$. Summation over repeated indices is assumed.  Using for instance a completeness property of the Pauli matrices 
one can verify that
\begin{equation}
\sum_{i=1}^3 (x^i)^2 = 1 \label{extrinsic-2-sphere}
\end{equation}
For further reference, let us also record here two other direct consequences of equation (\ref{stereo}), namely
\begin{equation}
\frac{\partial x^i}{\partial \pi^A} = \frac{\pi_A\bar{\pi}_B\bar{\pi}_D(\sigma^i)^{BD}}{(\bar{\pi}_C\pi^C)^2} \label{gjhjk}
\end{equation}
and 
\begin{equation}
\epsilon_{ijk\,}x^j\frac{\partial x^k}{\partial\pi^A} = i \frac{\partial x^i}{\partial\pi^A} \label{gdfgh}
\end{equation}
We conventionally lower the $SL(2,\mathbb{C})$ indices by means of the two-dimensional $\epsilon$-symbol and raise them by means of its inverse. For example, by definition, $\pi_A=\pi^B\epsilon_{BA}$  and so by way of consequence, $\pi^A = \epsilon^{AB}\pi_B$. Note incidentally that $(\sigma^i)^{AB} = (\sigma^i)^{BA}$ for all $i$. 

So, to each point of homogeneous coordinates $\pi^A$ on the Riemann sphere we can associate in $\mbox{End}\,T({\cal M})$ the complex structure 
\begin{equation}
J(\pi,\bar{\pi}) = x^i J_i \label{jay-pi-pibar}
\end{equation}
with $x^i$ given by (\ref{stereo}). By resorting to the quaternionic relations (\ref{quater-neon}) as well as to the equations (\ref{extrinsic-2-sphere}) and (\ref{gdfgh}) one can show that 
\begin{equation}
[I + i J(\pi,\bar{\pi})]\frac{\partial J(\pi,\bar{\pi})}{\partial\pi^A} = \frac{\partial J(\pi,\bar{\pi})}{\partial\pi^A}[I - i J(\pi,\bar{\pi})] = 0 \label{proj-dJ}
\end{equation}

The twistor space $Z$ is defined to be the direct product manifold ${\cal M} \times \mathbb{CP}^1$ endowed with the following complex structure
\begin{equation}
{\cal J} = J(\pi,\bar{\pi}) + i\frac{\partial}{\partial\pi^A} \otimes d\pi^A - i \frac{\partial}{\partial\bar{\pi}_A}\otimes d\bar{\pi}_A \label{jay-zee}
\end{equation}
This is an element of $\mbox{End}\,[T({\cal M})\oplus T(\mathbb{CP}^1)] \simeq \mbox{End}\,T({\cal M}) \oplus \mbox{End}\,T(\mathbb{CP}^1)$, with the first component given by the complex structure (\ref{jay-pi-pibar}) and the second by the standard complex structure on $\mathbb{C}^2$ that descends on $\mathbb{CP}^1$. The twistor space can be viewed as a holomorphic fibration over $\mathbb{CP}^1$, with the fiber corresponding to the point of homogeneous coordinates $\pi^A$ being  a copy of the manifold ${\cal M}$ endowed with the complex structure $J(\pi,\bar{\pi})$. The holomorphic sections of this fibration are termed twistor lines. 

To prove that this is indeed a complex structure, consider two arbitrary vector fields from $T(Z) \simeq T({\cal M})\oplus T(\mathbb{CP}^1)$
\begin{equation}
{\cal X} = X + X^A\frac{\partial}{\partial\pi^A} + \bar{X}_A\frac{\partial}{\partial\bar{\pi}_A} \qquad\mbox{and}\qquad
{\cal Y} = Y \,+ Y^A\frac{\partial}{\partial\pi^A} + \,\bar{Y}_A\frac{\partial}{\partial\bar{\pi}_A} 
\end{equation}
where $X, Y \in T({\cal M})$.  A direct calculation shows that the Nijenhuis tensor on $Z$ corresponding to the almost complex structure (\ref{jay-zee}) takes, when evaluated on ${\cal X}$ and ${\cal Y}$, the following form
\begin{eqnarray}
N_{{\cal J}}({\cal X},{\cal Y}) &\hspace{-7pt} = \hspace{-7pt}& N_{J(\pi,\bar{\pi})}(X,Y) \nonumber \\ [4pt] 
&\hspace{-7pt} - \hspace{-7pt}& i X^A[I+iJ(\pi,\bar{\pi})] \frac{\partial J(\pi,\bar{\pi})}{\partial\pi^A}Y + i Y^A[I+iJ(\pi,\bar{\pi})] \frac{\partial J(\pi,\bar{\pi})}{\partial\pi^A}X \nonumber \\ [2pt]
&\hspace{-7pt} + \hspace{-7pt}&  i \bar{X}_A\,[I-iJ(\pi,\bar{\pi})] \frac{\partial J(\pi,\bar{\pi})}{\partial\bar{\pi}_A}Y -  i \bar{Y}_A\,[I-iJ(\pi,\bar{\pi})] \frac{\partial J(\pi,\bar{\pi})}{\partial\bar{\pi}_A}X \label{Nijenhuis-jay-zee}
\end{eqnarray}
The integrability of $J(\pi,\bar{\pi})$ on ${\cal M}$ together with equation (\ref{proj-dJ}) imply then that it vanishes. By the Newlander-Nirenberg theorem it follows that the complex structure (\ref{jay-zee}) is indeed integrable and that $Z$ is a complex manifold. As we shall see, the holomorphic structure on $Z$ turns out to encode all the metric information of the hyperk\"ahler manifold.
  
 Let $\omega_1$, $\omega_2$, $\omega_3$ be the K\"ahler 2-forms corresponding to the complex structures $J_1$, $J_2$, $J_3$ on ${\cal M}$ and define the 2-form-valued ${\cal O}(2)$ section 
\begin{equation}
\omega(\pi) = -\frac{1}{2}\pi^A\pi^B(\sigma^i)_{AB}\,\omega_i \label{twisted-hol-2-form}
\end{equation}
Let also $g$ be the metric on ${\cal M}$. Then, for any two vector fields $X, Y \in T({\cal M})$, we have
 \begin{equation}
\omega(\pi)(X,Y) = -\frac{1}{2}\pi^A\pi^B(\sigma^i)_{AB}\,g(X,J_iY) = \frac{1}{2}\bar{\pi}_C\pi^C\pi_A\,g(X,\frac{\partial J(\pi,\bar{\pi})}{\partial \bar{\pi}_A}Y)
\end{equation}
The second equality follows from the equation (\ref{gjhjk}). Based on this result and the hermiticity of the metric we may write
\begin{equation}
\omega(\pi)([I+iJ(\pi,\bar{\pi})]X,Y) = \frac{1}{2}\bar{\pi}_C\pi^C\pi_A\,g(X,[I-iJ(\pi,\bar{\pi})]\frac{\partial J(\pi,\bar{\pi})}{\partial \bar{\pi}_A}Y) = 0
\end{equation}
The last equality follows from the equation (\ref{proj-dJ}).  We have thus shown that $\omega(\pi)$ is a $(2,0)$-form on ${\cal M}$ with respect to the complex structure $J(\pi,\bar{\pi})$. Since all three K\"ahler 2-forms $\omega_i$ are closed, it follows that  $\omega(\pi)$ is also closed and, moreover, holomorphic with respect to $J(\pi,\bar{\pi})$. In other words, for each point labelled by $\pi=(\pi^1,\pi^2)$ in $\mathbb{CP}^1$, $\omega(\pi)$ is a holomorphic symplectic form on the fiber of the projection $Z \rightarrow \mathbb{CP}^1$ above $\pi$.
\begin{figure}[htb]
\centering
\scalebox{0.50}{\includegraphics{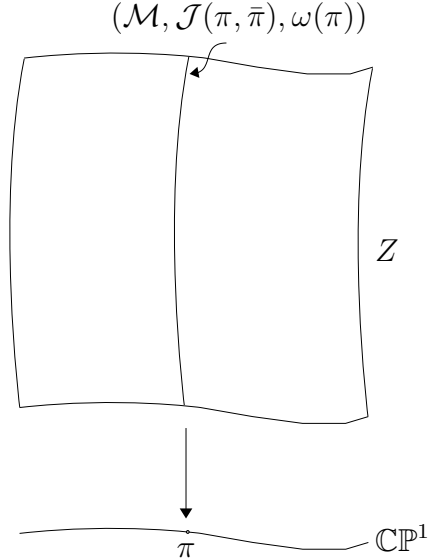}}
\put(-100,199){$({\cal M},{\cal J}(\pi,\bar{\pi}),\omega(\pi))$}
\put(0,110){$Z$}
\put(0,0){$\mathbb{CP}^1$}
\put(-75,-2){$\pi$}
\caption{The twistor fibration}
\end{figure} 

The antipodal conjugation on the sphere of hyperk\"ahler complex structures induces on the twistor space $Z$ the $\mathbb{Z}_2$-action
\begin{equation}
(m,\pi^A) \stackrel{\rm a.\,c.}{\longrightarrow} (m,\bar{\pi}^A)  \label{a-podal-conj}
\end{equation}
for any $(m,\pi^A) \in {\cal M}\times\mathbb{CP}^1$. This action maps ${\cal J}$ into $-{\cal J}$ and therefore defines a {\it real structure} on $Z$. In particular, note that (\ref{a-podal-conj}) takes $\omega(\pi)$ into its complex conjugate. This can be seen immediately by resorting to the hermiticity property $\overline{(\sigma^i)}_{AB} = (\sigma^i)^{BA}$.

\section{The generalized Legendre transform}

In the previous section we took a global approach and used homogeneous coordinates to describe the $\mathbb{CP}^1$ base of the twistor fibration. It is beneficial to take also a local point of view and describe bundles over $\mathbb{CP}^1$ and their sections in terms of patching by means of transition functions. Let $U_1$ and $U_2$ be the standard open charts on $\mathbb{CP}^1$, corresponding to the inhomogeneous coordinates $\zeta = \pi^2/\pi^1$ and $\tilde{\zeta} = \pi^1/\pi^2$, respectively. On $U_1\cap U_2$ these are related by $\tilde{\zeta} = 1/\zeta$. Choosing the trivialization $\pi = (1,\zeta)$, the 2-form-valued ${\cal O}(2)$ section (\ref{twisted-hol-2-form}) takes for example on $U_1$ the form
\begin{equation}
\omega(\zeta) = \omega_+ - \omega_3\zeta - \omega_-\zeta^2 \label{omega(zeta)}
\end{equation}
where $\omega_{\pm} = (\omega_1 \pm i\omega_2)/2$ are the hyperk\"ahler holomorphic and anti-holomorphic 2-forms, respectively. The reality condition induced by antipodal conjugation reads in local coordinates
\begin{equation}
\overline{\omega\!\left(-\frac{1}{\bar{\zeta}}\right)} = - \!\left(\frac{1}{\zeta}\right)^2\! \omega(\zeta)
\end{equation}
Consider now the fiber of $Z \rightarrow \mathbb{CP}^1$ above $\zeta \in U_1$. As we have seen, $\omega(\zeta)$ is a holomorphic symplectic form with respect to the distinguished complex structure on the fiber. In particular, this means that one can always introduce a set of complex holomorphic Darboux coordinates $p = p(\zeta)$ and $q = q(\zeta)$, in terms of which $\omega(\zeta) = dq \wedge dp$. If the real dimension of the hyperk\"ahler manifold is $4n$ then there are in fact $n$ such pairs of Darboux coordinates and the symplectic form is a direct sum of such terms. To decongest the notation we omit here the indices differentiating among these, but they should be understood to exist. Similarly, one can introduce complex holomorphic Darboux coordinates $P = P(\tilde{\zeta})$ and $Q = Q(\tilde{\zeta})$ on the fiber above $\tilde{\zeta} \in U_2$. If these can furthermore be chosen such that
\begin{equation}
Q(\tilde{\zeta}) = \overline{q\!\left(-\frac{1}{\bar{\zeta}}\right)} \qquad\mbox{and}\qquad
P(\tilde{\zeta}) = -\overline{p\!\left(-\frac{1}{\bar{\zeta}}\right)} \label{PQ-reality}
\end{equation}
then the required reality property of $\omega$  with respect to antipodal conjugation holds automatically. Since $\omega$ is ${\cal O}(2)$-valued, on $U_1 \cap U_2$ we may write 
\begin{equation}
\omega = dq \wedge dp = \zeta^2 dQ \wedge dP   \label{ttsm}
\end{equation}
The transition between the two sets of Darboux coordinates can be thus viewed as a twisted holomorphic symplectomorphism.

At this point we make a few restrictive assumptions. First of all, we assume that the $q$-coordinates are ${\cal O}(2j)$ sections over $\mathbb{CP}^1$ (not necessarily with the same value of $j$, when there are several of these) and that the $p$-coordinates are non-singular close to $\zeta = 0$. Accordingly, we have the following power expansions
\begin{equation}
q = \sum_{n=0}^{2j} q_n\zeta^{n}  \qquad\mbox{and}\qquad
p = \sum_{n=0}^{\infty} p_n\zeta^n \label{pq-pwr-exps}
\end{equation}
The twisted holomorphic symplectomorphism (\ref{ttsm}) can be derived from the twisted type II canonical transformation 
\begin{eqnarray}
Q &\hspace{-7pt} = \hspace{-7pt}& \,\frac{(-)^j}{\zeta^{2j}}\,\frac{\partial F_2(q,P)}{\partial P} \label{typeII-rel1} \\ [2pt]
p &\hspace{-7pt} = \hspace{-7pt}& \frac{(-)^j}{\zeta^{2j-2}}\frac{\partial F_2(q,P)}{\partial q} \label{typeII-rel2}
\end{eqnarray}
Our second assumption is that the generating function is of the form
\begin{equation}
F_2(q,P) = qP + (-)^j\zeta^{2j} G(q)
\end{equation}
with $G$ depending solely on $q$ and perhaps on $\zeta$ but not on $P$. The factor in front of $G$ is nonessential, it serves a formal purpose and it can be in principle absorbed in the definition of $G$. Equation (\ref{typeII-rel1}) encodes in fact  a reality condition for the ${\cal O}(2j)$ sections. Indeed, together with the first equation (\ref{PQ-reality}) it yields that
\begin{equation}
\overline{q\!\left(-\frac{1}{\bar{\zeta}}\right)}  = (-)^j\! \left(\frac{1}{\zeta}\right)^{2j} \!q(\zeta)
\end{equation}
In terms of the parameters, this is equivalent to $\bar{q}_n = (-)^{j-n} q_{2j-n}$. In particular, the middle parameter $q_j$ is always real. Equation (\ref{typeII-rel2}), on the other hand, yields a patching formula for the $p$-coordinates on $U_1 \cap U_2$. Integrating it on a closed contour around $\zeta = 0$ we obtain
\begin{equation}
\frac{\partial F}{\partial q_n} = p_{1-n}-(-)^{j-n}\bar{p}_{1-2j+n} \label{dFdqm}
\end{equation}
for $n=1,\cdots,2j$, where we have defined
\begin{equation}
F = \frac{1}{2\pi i}\oint d\zeta\, G(q(\zeta)) \label{F-oint-G}
\end{equation}
The integration contour should be chosen such that the resulting $F$ be purely imaginary. Explicitly, the system of equations (\ref{dFdqm}) reads
\begin{eqnarray}
\frac{\partial F}{\partial q_0} &\hspace{-7pt} = \hspace{-7pt}& p_1 \label{L-rel-0} \\ [1pt]
\frac{\partial F}{\partial q_1} &\hspace{-7pt} = \hspace{-7pt}& p_0 + (-)^j\bar{p}_{2-2j} \label{L-rel-1} \\
\frac{\partial F}{\partial q_2} &\hspace{-7pt} = \hspace{-7pt}& \cdots = \frac{\partial F}{\partial q_j} = 0 \label{L-rel-2}
\end{eqnarray}
The remaining equations  are just the complex conjugates of these. Note that the second term in the r.h.s. of (\ref{L-rel-1}) vanishes for all values of $j$ except for $j=1$, when it is equal to $-\bar{p}_0$.

Substituting the power expansions (\ref{pq-pwr-exps}) into the Darbox form of $\omega$ on $U_1$, we get
\begin{equation}
\omega(\zeta) = dq \wedge dp = dq_0\wedge dp_0 + (dq_0\wedge dp_1 + dq_1\wedge dp_0)\zeta + \cdots
\end{equation}
Comparison with (\ref{omega(zeta)}) yields
\begin{eqnarray}
\omega_+ &\hspace{-7pt} = \hspace{-7pt}& dq_0\wedge dp_0 \\ [5pt]
\omega_3\, &\hspace{-7pt} = \hspace{-7pt}& -(dq_0\wedge dp_1 + dq_1\wedge dp_0)
\end{eqnarray}
Observe now that in the limit when $\zeta \rightarrow 0$, the complex structure $J(\zeta,\bar{\zeta})$ on the twistor fiber above $\zeta$ converges to $J_3$, while the holomorphic symplectic structure $\omega(\zeta)$ converges to $\omega_+$. Since $p$ and $q$ are holomorphic coordinates with respect to $J(\zeta,\bar{\zeta})$ we conclude that $q_0$ and $p_0$ are holomorphic with respect to $J_3$. Using that, we can write $\omega_3$ locally as a total derivative
\begin{equation}
\omega_3 = d(p_1dq_0 - q_1dp_0) \label{kask}
\end{equation}
On another hand, defining a real function $K$ by the Legendre transform
\begin{equation}
iK(q_0,\bar{q}_0,p_0,\bar{p}_0) = F(q_0,\bar{q}_0,q_1,\bar{q}_1,\cdots) - (p_0q_1-\bar{p}_0\bar{q}_1) \label{K-pot-L-tr}
\end{equation}
we get, by means of the equations (\ref{L-rel-0}) through (\ref{L-rel-2}), that
\begin{equation}
i \partial K = i\left(\frac{\partial K}{\partial q_0}dq_0 + \frac{\partial K}{\partial p_0}dp_0 \right) = p_1dq_0 - q_1dp_0
\end{equation}
Plugging this in (\ref{kask}) yields 
\begin{equation}
\omega_3 = i\bar{\partial}\partial K
\end{equation}
which means that the function $K$ thus defined is in fact a K\"ahler potential for the complex structure $J_3$.

\section{Hyperk\"ahler metrics from the generalized Legendre transform} \label{SEC:hk-metrics-GLT}

The equations (\ref{F-oint-G}) and (\ref{K-pot-L-tr}) together with (\ref{L-rel-1}) and (\ref{L-rel-2}) can be used to construct hyperk\"ahler metrics from holomorphic objects. We summarize them here, reformulating them at the same time for later convenience, by introducing new notations and performing several redefinitions. Thus, for instance, we will henceforth cast real ${\cal O}(2j)$ sections over $Z$ (we can define ${\cal O}(k)$ sections over $Z$ by pulling back from $\mathbb{CP}^1$) or, equivalently, ${\cal O}(2j)$ multiplets, not in the usual polynomial form but in the following "tropical" form
\begin{equation}
\eta^{(2j)}(\zeta)=\frac{\bar{z}}{\zeta^j} + \frac{\bar{v}}{\zeta^{j-1}} + \frac{\bar{t}}{\zeta^{j-2}} + \cdots + x + (-)^j(\cdots + t \zeta^{j-2} - v\zeta^{j-1} + z \zeta^j) 
\end{equation}
The reality condition with respect to antipodal conjugation reads now
\begin{equation}
\eta^{(2j)}(-\frac{1}{\bar{\zeta}}) = \overline{\eta^{(2j)}(\zeta)}
\end{equation}
For simplicity, the $G$-functions that we will consider in this section will apparently depend on only one ${\cal O}(2j)$ multiplet. Nevertheless, the subsequent discussion can be straightforwardly generalized to include several multiplets or combinations of multiplets with different values of $j$, including $j=1$. This will amount, essentially, to introducing additional indices to differentiate among these. If we define
\begin{equation}
F = \oint \frac{d\zeta}{\zeta} G(\eta^{(2j)}(\zeta)) \label{F-pot}
\end{equation}
choosing the contour such that the outcome of the integration is real (take note of the difference with respect to (\ref{F-oint-G})!) then the generalized Legendre transform is given by
\begin{equation}
K(z,\bar{z},u,\bar{u})=F(z,\bar{z},v,\bar{v},t,\bar{t}, \cdots \!, x) - uv- \bar{u}\bar{v} \label{K-pot}
\end{equation}
together with the conditions
\begin{eqnarray}
\frac{\partial F}{\partial v} &\hspace{-7pt} = \hspace{-7pt}& u \label{Legendre-rel1} \\ [4pt]
\frac{\partial F}{\partial t}  &\hspace{-7pt} = \hspace{-7pt}& \cdots  = \frac{\partial F}{\partial x} = 0 \label{Legendre-rel2}
\end{eqnarray}
If the multiplet on which the function $G$ depends is of ${\cal O}(2)$ type, then $v=\bar{v}=x$ and these equations are replaced by the single equation
\begin{equation}
\frac{\partial F}{\partial x} = u + \bar{u} \hspace{42pt}  \label{Legendre-rel-O2}
\end{equation}
The generalized Legendre transform yields a K\"ahler potential $K$ as well as a corresponding system of holomorphic coordinates, $z$ and $u$. Equations (\ref{Legendre-rel1})-(\ref{Legendre-rel2}) can in principle be solved to give $v$, $\bar{v}$, $t$, $\bar{t}$, \dots, $x$ as implicit functions of $z$, $\bar{z}$, $u$, $\bar{u}$. Implicit differentiation returns
\begin{equation}
\frac{\partial a}{\partial z} = -F^{ab}F_{bz} \qquad\qquad
\frac{\partial a}{\partial u} = F^{av}  \label{dadz-dadu}
\end{equation}
where $a$ and $b$ run over the values $v$, $\bar{v}$, $t$, $\bar{t}$, \dots, $x$ and $F^{ab}$ is by definition the inverse of the matrix of second derivatives $F_{ab}$. Summation over repeated indices is assumed. On another hand, taking the derivatives of (\ref{K-pot}) with respect to the holomorphic coordinates and imposing afterwards the generalized Legendre relations (\ref{Legendre-rel1})-(\ref{Legendre-rel2}), one gets that
\begin{equation}
\frac{\partial K}{\partial z} = \frac{\partial F}{\partial z} \qquad\mbox{and}\qquad \frac{\partial K}{\partial u} = -v \label{Kz-Ku}
\end{equation}
Using further equations (\ref{dadz-dadu}) to take the derivatives of (\ref{Kz-Ku}) with respect to the anti-holomorphic variables this time, one obtains the metric components 
\begin{eqnarray}
\left( 
\begin{array}{cc}
K_{z\bar{z}} & K_{z\bar{u}} \\
K_{u\bar{z}} & K_{u\bar{u}} \\
\end{array}
\right)
&\hspace{-7pt} = \hspace{-7pt}&
\left( 
\begin{array}{ll}
F_{z\bar{z}}-F_{za}F^{ab}F_{b\bar{z}} & F_{za}F^{a\bar{v}}  \\
F^{va}F_{a\bar{z}} & -F^{v\bar{v}} \\
\end{array}
\right) \label{metric-comps1}
\end{eqnarray}
In the generalized Legendre transform approach the three standard K\"{a}hler forms are given by
\begin{equation}
\omega^3=-i\,\partial\bar{\partial}K \qquad\qquad \omega^+=dz\wedge du \qquad\qquad \omega^-=\overline{\omega^+\rule{-1.3mm}{0mm}} \label{HK-forms}
\end{equation}
From the last two expressions in (\ref{HK-forms}) and the fact that for any hyperk\"{a}hler manifold, in a  coordinate basis holomorphic with respect to the complex structure $J_3$, the components of its (2,0) and (0,2) K\"{a}hler forms satisfy
\begin{equation}
\omega^+{}_{\mu\rho} \, \omega^{-\rho\nu} = - \delta_\mu{}^\nu \qquad\qquad
\omega^{-\rho\sigma} K_{\rho\bar{\mu}}  K_{\sigma\bar{\nu}} = \omega^-{}_{\bar{\mu}\bar{\nu}} 
\end{equation}
it follows that the inverse metric relates in a direct way the metric itself, namely
\begin{equation}
\left( 
\begin{array}{cc}
K^{z\bar{z}} & K^{z\bar{u}} \\
K^{u\bar{z}} & K^{u\bar{u}} \\
\end{array}
\right) =
\left( \hspace{-4pt}
\begin{array}{rr}
K_{u\bar{u}} \!\!&\!\! - K_{u\bar{z}} \\
- K_{z\bar{u}} \!\!&\!\! K_{z\bar{z}} \\
\end{array}
\!\right) \label{inv-metric}
\end{equation}
On the other hand, one can attempt to invert directly (\ref{metric-comps1}) using the following elementary linear algebra result: if $A$ and $D$ are non-singular square matrices then
\begin{equation}
\left(
\begin{array}{cc}
A & B \\ 
C & D
\end{array}
\right)^{\!\!-1}
=
\left(\!\!\!
\begin{array}{cc}
(A-BD^{-1}C)^{-1} 			\!\!\!&\!\! -A^{-1}B(D-CA^{-1}B)^{-1} \\ 
-D^{-1}C(A-BD^{-1}C)^{-1} 	\!\!\!&\!\! (D-CA^{-1}B)^{-1}
\end{array}
\!\!\right)
\end{equation}
By identifying the upper-left block of the inverse metric obtained in this way to the upper-left block of (\ref{inv-metric}), one gets
\begin{equation}
F_{z\bar{z}} - F_{za}F^{ab}F_{b\bar{z}} + F_{za}F^{a\bar{v}}H_{\bar{v}v}F^{vb}F_{b\bar{z}} = - H_{\bar{v}v} 
\end{equation}
where $H_{\bar{v}v}$ denotes the matrix inverse  of $F^{v\bar{v}}$. This allows then one to rewrite (\ref{metric-comps1}) in the more symmetric form
\begin{equation}
\left( 
\begin{array}{cc}
K_{z\bar{z}} & K_{z\bar{u}} \\
K_{u\bar{z}} & K_{u\bar{u}} 
\end{array}
\right)
=
\left(\! 
\begin{array}{cc}
\ \mathbb{I}  & \ F_{za}F^{a\bar{v}} \\
\ 0  & \ -F^{v\bar{v}} 
\end{array}
\!\right)
\!\left(\!\! 
\begin{array}{cc}
-H_{\bar{v}v} \!&\! 0 \\
0 \!&\! -H_{\bar{v}v} 
\end{array}
\!\!\right)
\!\left(\! 
\begin{array}{cc}
\mathbb{I} & 0 \\ 
F^{va}F_{a\bar{z}} \!&\! -F^{v\bar{v}} 
\end{array}
\!\right) 
\label{metric-comps2}
\end{equation}
or, equivalently,
\begin{equation}
ds^2 = - dzH_{\bar{v}v}d\bar{z} - (du - dz F_{za}F^{a\bar{v}}H_{\bar{v}v})F^{v\bar{v}}(d\bar{u}-H_{\bar{v}v}F^{va}F_{a\bar{z}}d\bar{z}) \label{GLT-metric}
\end{equation}
From the form (\ref{metric-comps2}) of the metric it is immediately apparent that the following Monge-Amp\`{e}re equation holds
\begin{equation}
\det K_{(z,u)} = 1
\end{equation}
This result was established by Lindstr\"{o}m and Ro\v{c}ek \cite{Lindstrom:1983rt} for the case of several ${\cal O}(2)$ multiplets and subsequently by Cherkis \cite{Cherkis-unpb} for the case of one ${\cal O}(2j)$ multiplet. Both proofs exploit the fact that the contour-integral form (\ref{F-pot}) implies that the function $F$ satisfies the following set of second order differential equations
\begin{eqnarray}
F_{z\bar{z}} &\hspace{-7pt} = \hspace{-7pt}& - F_{v\bar{v}} = F_{t\bar{t}\,} = \cdots = (-)^j F_{xx} \nonumber \\ [4pt]
F_{z\bar{v}} &\hspace{-7pt} = \hspace{-7pt}& - F_{v\bar{t}\,} = \cdots \nonumber \\ [4pt]
F_{zt\,}  &\hspace{-7pt} = \hspace{-7pt}& F_{vv} \hspace{24.5pt} \mbox{etc.} \nonumber \\ [4pt]
F_{zv} &\hspace{-7pt} = \hspace{-7pt}& F_{vz} \hspace{25pt} \mbox{etc.} \label{diff-eqs}
\end{eqnarray}
The argument presented here can be adapted in a straightforward manner to generic combinations of multiplets. 

In the case of a single ${\cal O}(2j)$ multiplet the components of the metric can be expressed as ratios of determinants of Hankel matrices, as follows from
\begin{eqnarray}
F_{za}F^{a\bar{v}} &\hspace{-7pt} = \hspace{-7pt}& \frac{(-)^j}{\det F} \det\, (h_{a+b-1})_{-(j-1) \leq a,b \leq +(j-1)}  \\
-F^{v\bar{v}} &\hspace{-7pt} = \hspace{-7pt}& \frac{(-)^j}{\det F} \det\, (h_{a+b-1})_{-(j-2) \leq a,b \leq +(j-1)} \\ [6pt]
\det F &\hspace{-7pt} = \hspace{-7pt}& \det\, (h_{a+b})_{-(j-1) \leq a,b \leq +(j-1)} 
\end{eqnarray}
where
\begin{equation}
h_k = \oint \frac{d\zeta}{\zeta} \zeta^{-k} \frac{\partial^2G}{\partial\eta^2}
\end{equation}
These relations can be derived by writing the elements of the matrix $F^{ab}$ that occur in the l.h.s. in terms of the cofactors of its inverse matrix, $F_{ab}$, and,  in the first case only, by using subsequently Laplace's determinant expansion formula. The resulting determinants can then be algebraically manipulated into the forms displayed above. This generalizes an observation made by Cherkis and Hitchin in \cite{Cherkis:2003wk}.

\section{Pure ${\cal O}(2)$ constructions} \label{SEC:O2-GLT}

\subsection{The generic case}

As noted already, ${\cal O}(2)$ multiplets are in some sense special and their presence require certain adjustments to the generalized Legendre transform formalism developed above. There are no extremization conditions of the type (\ref{Legendre-rel2}) associated with ${\cal O}(2)$ multiplets as they do not have enough components, or,  in a field-theoretical language, there are no auxiliary fields to be integrated out. There is only one Legendre relation associated to each of them. Distinguishing with an index $I$ between the various ${\cal O}(2)$ multiplets of the theory, these take the form
\begin{equation}
\frac{\partial F}{\partial x^I} = u_I + \bar{u}_I \label{L-rel-O2}
\end{equation}
The departure of these relations from the regular form (\ref{Legendre-rel1}) for  higher-degree multiplets stems from  the reality properties of ${\cal O}(2)$ multiplets which require that $x^I$ be real. In the absence of $\bar{u}_I$ the regular form would simply be inconsistent. 

Let us consider the class of  $F$-functions constructed exclusively out of a number of $n$ ${\cal O}(2)$ multiplets and review the particularities of the corresponding $4n$-real dimensional hyperk\"{a}hler manifolds.  This class of manifolds was studied in \cite{Hitchin:1986ea} and \cite{MR953820}. The metric formula (\ref{GLT-metric}) specializes in this case to
\begin{equation}
ds^2 = -dz^IF_{x^Ix^J}d\bar{z}^J - (du_I-dz^KF_{z^Kx^I})F^{x^Ix^J}(d\bar{u}_J-F_{x^J\bar{z}^L}d\bar{z}^L) \label{GLT-metric-O2}
\end{equation}
where $F^{x^Ix^J}$ is the matrix inverse of $F_{x^Ix^J}$, whereas the holomorphic $(2,0)$-form from (\ref{HK-forms}) becomes
\begin{equation}
\omega^+ = dz^I\wedge du_I \label{hol-2-form-O2}
\end{equation}
The fact that the holomorphic coordinates $u_I$ occur exclusively in the combination $u_I+\bar{u}_I$ implies that the metric is independent of the imaginary part of $u_I$, in other words it has  $n$ abelian holomorphic Killing vectors
\begin{equation}
\tilde{X}^I = i\left(\frac{\partial}{\partial u_I}-\frac{\partial}{\partial\bar{u}_I}\right) \label{O2-3-hol-isometry}
\end{equation}
Clearly, these Killing vectors preserve the $(2,0)$-form (\ref{hol-2-form-O2}) and hence they are not only holomorphic but in fact tri-holomorphic.

In the holomorphic coordinate basis $z^I$, $u_I$ the hyperk\"{a}ler structure is manifest but the underlying $SO(3)$ structure deriving from the ${\cal O}(2)$ multiplets is obscure and so are the abelian  isometries. We can make the symmetries manifest and obscure the holomorphic structure by switching to a set of real coordinates defined as follows
\begin{eqnarray}
\vec{r}^{\,I} &\hspace{-7pt} = \hspace{-7pt}& (z^I+\bar{z}^I, -i(z^I-\bar{z}^I),\, x^I) \label{O2-radii} \\ [3pt]
\psi_I &\hspace{-7pt} = \hspace{-7pt}& \mbox{Im}\, u_I
\end{eqnarray}
In this coordinate basis the metric takes a generalized Gibbons-Hawking form,
\begin{equation}
ds^2 \sim \Phi_{IJ}\, d\vec{r}^{\,I} \!\cdot\! d\vec{r}^{\,J} + (\Phi^{-1})^{IJ}(d\psi_I+\vec{A}_{IK} \!\cdot\! d\vec{r}^{\,K})(d\psi_J+\vec{A}_{JL} \!\cdot\! d\vec{r}^{\,L}) \label{Gibbons-Hawking}
\end{equation}
where the tilde symbolizes "equal, up to an overall factor $1/2$" and $\vec{A}$, $\Phi$ are defined by
\begin{eqnarray}
\vec{A}_{IJ} \!\cdot\! d\vec{r}^{\,J} &\hspace{-7pt} = \hspace{-7pt}& \frac{i}{2}(F_{x^Iz^J}dz^J - F_{x^I\bar{z}^J}d\bar{z}^J) \label{A_IJ} \\ [2pt]
\Phi_{IJ} &\hspace{-7pt} = \hspace{-7pt}& -\frac{1}{2}F_{x^Ix^J} \label{phi_IJ}
\end{eqnarray}
To derive this expression of the metric from the form (\ref{GLT-metric-O2}) it is necessary to use the fact that
\begin{equation}
F_{x^Iz^J} = F_{x^Jz^I}
\end{equation}
a direct consequence of the contour-integral form (\ref{F-pot}) of $F$. By resorting again to this relation, this time in conjunction with equations (\ref{A_IJ}) and (\ref{phi_IJ}), one can prove moreover that the following generalized Bogomol'nyi equations hold
\begin{equation}
\vec{\nabla}_I \times \vec{A}_{KJ} = - \vec{\nabla}_I \Phi_{KJ} \qquad\mbox{and}\qquad \vec{\nabla}_I \Phi_{KJ} = \vec{\nabla}_J \Phi_{KI}
\end{equation}
where $\vec{\nabla}_I = \partial/\partial \vec{r}^{\,I}$ are $\mathbb{R}^3$ gradient operators.

\subsection{The Taub-NUT metric}

As an example, we review the construction of the 4-dimensional hyperk\"ahler manifold generated by the $F$-function
\begin{equation}
F = - \frac{1}{2\pi i h}\oint_{\Gamma_0} \frac{d\zeta}{\zeta} (\eta^{(2)})^2 + \oint_{\Gamma} \frac{d\zeta}{\zeta} \eta^{(2)} \ln \eta^{(2)}
\end{equation}
This gives the Taub-NUT metric when $h$ is finite and the flat space metric when $h \rightarrow \infty$ \cite{Karlhede:1984vr}. The generalized Legendre relation (\ref{L-rel-O2}) becomes in this case
\begin{equation}
\frac{x}{h} + \tanh^{-1} \frac{x}{r} = -\frac{u+\bar{u}}{2}
\end{equation}
This equation is transcendental for finite $h$ and algebraic for $h \rightarrow \infty$. It can be used to express the real parameter $x$ as a function (explicitly in the latter case, implicitly in the former) of the holomorphic coordinates $z$, $u$ and their complex conjugates. The Legendre-Fenchel transform yields the K\"ahler potential 
\begin{eqnarray}
K = 2r + \frac{r^2-2|z|^2}{h} 
\end{eqnarray}
In a holomorphic coordinate basis the metric takes the form
\begin{equation}
ds^2 = 2\Phi\, dz d\bar{z} + (2\Phi)^{-1} (du - A_z dz)(d\bar{u} - A_{\bar{z}}d\bar{z}) \label{A-metric-hol}
\end{equation}
with
\begin{equation}
\Phi = \frac{1}{h} + \frac{1}{r} \qquad\mbox{and}\qquad A_{z} = \frac{x}{zr}
\end{equation}
Besides the $U(1)$ tri-holomorphic isometry (\ref{O2-3-hol-isometry}), this metric possesses an additional non-tri-holomorphic  $SO(3)$ isometry induced by the natural $SO(3)$ action on the parameter space of the ${\cal O}(2)$ real section, generated by the vector fields
\begin{eqnarray}
X_-	&\hspace{-7pt} = \hspace{-7pt}& x \frac{\partial}{\partial z} + \frac{hr+2|z|^2}{hz}\frac{\partial}{\partial u} + 2 \frac{\bar{z}}{h}\frac{\partial}{\partial \bar{u}} \\ [1pt]
X_+ 	&\hspace{-7pt} = \hspace{-7pt}& x \frac{\partial}{\partial \bar{z}} + \frac{hr+2|z|^2}{h\bar{z}}\frac{\partial}{\partial \bar{u}} + 2\frac{z}{h} \frac{\partial}{\partial u} \\
X_{3}\, &\hspace{-7pt} = \hspace{-7pt}& -2i\left(z\frac{\partial}{\partial z} - \bar{z}\frac{\partial}{\partial \bar{z}}\right)
\end{eqnarray}
One can verify directly that these satisfy indeed the $SO(3)$ algebra
\begin{equation}
[X_i,X_j] = 2\epsilon_{ijk\,}X_k
\end{equation}
and that, moreover, they rotate the hyperk\"{a}ler 2-forms, {\it i.e.},
\begin{equation}
{\cal L}_{X_i} \omega_j = 2\epsilon_{ijk\,} \omega_k 
\end{equation}
as one would expect from general considerations. We use here the standard notation conventions $\displaystyle{X_{\pm} = (X_1\pm iX_2)/2}$ and $\displaystyle{\omega^{\pm} = (\omega_1\pm i\omega_2)/2}$. The actions of the isometry generators are hamiltonian with respect to their corresponding K\"{a}hler 2-forms, with moment maps
\begin{eqnarray}
\mu_1 &\hspace{-7pt} = \hspace{-7pt}& 2r + \frac{r^2-2|z|^2-2\mbox{Re}(z^2)}{h} \\ 
\mu_2 &\hspace{-7pt} = \hspace{-7pt}& 2r + \frac{r^2-2|z|^2+2\mbox{Re}(z^2)}{h} \\ 
\mu_3 &\hspace{-7pt} = \hspace{-7pt}& 2r + \frac{4|z|^2}{h}
\end{eqnarray}
If one writes the components of $\vec{r}$ in polar coordinates then the directional $\phi$ and $\theta$ angles associate in a natural manner with $\psi = \mbox{Im\,} u$ to parametrize the $SO(3)$ group manifold. One obtains in this way the more familiar form of the Taub-NUT metric
\begin{equation}
ds^2 \sim \Phi\, [dr^2+r^2(\sigma_1^2+\sigma_2^2)] + \Phi^{-1}\sigma_3^2
\end{equation}
with $\sigma_1$, $\sigma_2$, $\sigma_3$ representing the Cartan-Maurer left-invariant  $SO(3)$ 1-forms.

\section{The ${\cal O}(4)$ spectral curve} \label{SEC:O4-multiplets}

\subsection{Majorana normal form}

We will now focus our attention on the multiplets of ${\cal O}(4)$ type. These can be cast, generically, in either one of the following two local forms, to which we will henceforth refer as Majorana normal forms \cite{Majorana:1932ga}
\begin{eqnarray}
\eta^{(4)}(\zeta) &\hspace{-7pt} = \hspace{-7pt}& \frac{\bar{z}}{\zeta^2}+\frac{\bar{v}}{\zeta}+x-v\zeta+z\zeta^2 \nonumber \\
&\hspace{-7pt} = \hspace{-7pt}&  \frac{\rho}{\zeta^2} \frac{(\zeta-\alpha)(\bar{\alpha}\zeta+1)}{1+|\alpha|^2}\frac{(\zeta-\beta)(\bar{\beta}\zeta+1)}{1+|\beta|^2} \label{eta4}
\end{eqnarray}
The reality constraint satisfied by the multiplet is preserved only by the $PSU(2)$ subgroup of the $PSL(2,\mathbb{C})$ group of automorphisms of $\mathbb{CP}^1$ which acts on the inhomogeneous coordinate $\zeta$ by M\"obius transformations. The projective component of $PSU(2)$ generates a real scaling transformation.  Under the $SU(2)$ action, the polynomial coefficients transform together in the spin-$2$ irrep. On the other hand, the roots $\alpha$, $\beta$ and their antipodal conjugates can be viewed as points on a Riemann sphere endowed with the Fubini-Study metric.  The $SU(2)$ transformations correspond to isometric rotations of the sphere, under which the system of points moves together rigidly. For a more detailed discussion of the rotational properties of multiplets we refer the reader to the second article in this series, \cite{MM2}.

There are two irreducible spherical invariants that one can associate to $\eta^{(4)}$. These can be taken to be
\begin{eqnarray}
g_2 &\hspace{-7pt} = \hspace{-7pt}& 4 |z|^2 + |v|^2 + \frac{1}{3} x^2  \label{g2_Majorana} \\
g_3 &\hspace{-7pt} = \hspace{-7pt}&  \frac{8}{3} |z|^2 x - \frac{1}{3} |v|^2 x - \frac{2}{27} x^3 - z \bar{v}^2 - \bar{z} v^2 \label{g3_Majorana}
\end{eqnarray}
In principle, one can construct out of the Majorana coefficients other spherical invariants as well, but these turn out to be all reducible, in the sense that they can be decomposed into polynomial expressions in $g_2$ and $g_3$, with rational coefficients. If we substitute into these formulas the Vi\`ete expressions of $z$, $\bar{z}$, $v$, $\bar{v}$ and $x$ in terms of the roots and the scale factor, we get that
\begin{eqnarray}
g_2 &\hspace{-7pt} = \hspace{-7pt}& \ \frac{1}{3} \ \rho^2 (1 - k^2 + k^4) \label{r2}  \label{g2_rho} \\
g_3 &\hspace{-7pt} = \hspace{-7pt}& \frac{1}{27} \,\rho^3 (k^2 - 2) (2 k^2-1) (k^2+1) \label{g3_rho}
\end{eqnarray}
where 
\begin{equation}
k = \frac{|1+\bar{\alpha}\beta|}{\sqrt{(1+|\alpha|^2)(1+|\beta|^2)}} \label{k-modulus}
\end{equation}
is directly related to the Fubini-Study distance between $\alpha$ and $\beta$ (more precisely, the latter is given by $2\arccos k$).

The fact that $\eta^{(4)}$ has scaling weight 2 suggests considering the associated quartic plane curve
\begin{equation}
\eta^2 = \zeta^2\eta^{(4)}(\zeta) \label{Majorana_n.f.}
\end{equation}
The projection $(\zeta,\eta) \longmapsto \zeta$ is a two-sheeted branched covering of the Riemann sphere, with the holomorphic elliptic involution $(\zeta,\eta) \longmapsto (\zeta,-\eta)$ interchanging the two sheets except at its fixed points, which are branching points for the covering. When these are all different, the curve is non-singular. The curve has, additionally, two anti-holomorphic involutions or real structures $(\zeta,\eta) \longmapsto ( -1/\bar{\zeta},\pm\bar{\eta}/\bar{\zeta}^2)$ induced by the antipodal map on the sphere,  conjugated by the elliptic involution and preserving the set of branching points. This makes it a double-cover of the real projective plane, $\mathbb{R}\mathbb{P}^2 \simeq \mathbb{CP}^1/\mathbb{Z}_2$, and thus a real algebraic curve of genus 1. 

As an elliptic curve, it has an abelian differential form, {\it i.e.}, a globally defined holomorphic 1-form
\begin{equation}
\varpi = \frac{d\zeta}{2\zeta\sqrt{\eta^{(4)}(\zeta)}} \label{omega_Majorana}
\end{equation}

\subsection{Legendre normal form}

The equation of the curve can be cast into Legendre normal form by means of a birational transformation of $\zeta$ mapping three members of the root system $\alpha$,$-1/\bar{\alpha}$, $\beta$,$-1/\bar{\beta}$ to  $0$, $1$ and $\infty$, while also appropriately transforming $\eta$. As is well-known in the theory of elliptic curves, one can obtain in this way six possible moduli, representing points where the fourth root can be mapped. A very simple and elegant modulus, namely $k^2$, is returned by any one of the following four birational transformations: the one given in the form of the cross-ratio\footnote{The notation for cross-ratio that we employ here is $\displaystyle{[z_1,z_2,z_3,z_4] = \frac{(z_1-z_3)(z_2-z_4)}{(z_1-z_4)(z_2-z_3)}}$.}
\begin{equation} 
\nu = [\zeta,-\frac{1}{\bar{\alpha}},\ \alpha,\ \beta] \label{Majorana-Legendre}
\end{equation}
another one obtained by replacing in this relation $\alpha$ and $\beta$ with their antipodal conjugates, as well as two others that can be obtained from these two by interchanging $\alpha$ and $\beta$. Note, incidentally, that antipodal conjugation translates in this new context into complex-conjugation. In all these four cases the abelian 1-form (\ref{omega_Majorana}) transforms to
\begin{equation}
\varpi = \frac{d\nu}{2\sqrt{\rho\,\nu(\nu-1)(\nu-k^2)}} \label{omega_Legendre}
\end{equation}
and so equation~(\ref{Majorana_n.f.}) can be re-expressed in terms of $\nu$ and the coordinate $\mu = \eta\, \partial{\nu}/\partial{\zeta}$ as
\begin{equation}
\mu^2 = \rho\, \nu(\nu-1)(\nu-k^2) \label{Legendre_n.f.}
\end{equation}
The periods of $\varpi$ are obtained by integrating it over the canonical cycles. Through some standard changes of variables one can relate the period integrals to the familiar form of the complete Legendre elliptic integrals of the first kind. The period lattice reads
 \begin{equation}
\Lambda = \mathbb{Z}\cdot 2\omega + \mathbb{Z}\cdot 2\omega'
\end{equation}
with the half-periods given by
\begin{equation}
\omega = \frac{K(k)}{\sqrt{\rho}} \hspace{20pt} \mbox{and} \hspace{20pt} \omega' = \frac{iK(k')}{\sqrt{\rho}}  \label{half-periods}
\end{equation}
Since $0 < k, k' <1$, it follows that $K(k)$ and $K(k')$ are real and so $\Lambda$ is an orthogonal lattice. The Fubini-Study invariants $k$ and $k'$ play in this context the role of elliptic modulus and complementary modulus, respectively \cite{Ivanov:1995cy}.

\subsection{Weierstrass normal form}

At this point, let us note that equations~(\ref{g2_rho}) and (\ref{g3_rho}) can be re-written in terms of
\begin{equation}
e_1 = -\frac{\rho}{3} (k^2-2) \qquad\quad e_2 = \frac{\rho}{3} (2k^2-1) \qquad\quad e_3 = -\frac{\rho}{3}(k^2+1) \label{Weierstrass_roots}
\end{equation}
satisfying
\begin{equation}
e_1+e_2+e_3 = 0 \label{G1}
\end{equation}
as follows
\begin{eqnarray}
g_2 &\hspace{-7pt} = \hspace{-7pt}& -(e_1e_2+e_2e_3+e_3e_1) \label{G2} \\ [4pt]
g_3 &\hspace{-7pt} = \hspace{-7pt}& e_1e_2e_3 \label{G3}
\end{eqnarray}
Equations (\ref{Weierstrass_roots}) can be reverted to yield
\begin{equation}
k^2 = \frac{e_2-e_3}{e_1-e_3} \hspace{20pt} \mbox{and} \hspace{20pt} \rho = e_1-e_3
\end{equation}
Recognizing that the equations (\ref{G1}) through (\ref{G3}) are Vi\`{e}te-type formulas, it is then immediately apparent that $e_1$, $e_2$ and $e_3$ have to be the three roots of the cubic equation $X^3-g_2X-g_3 = 0$. This suggests that one should try to cast the equation of the elliptic curve associated to $\eta^{(4)}$ into the Weierstrass normal form
\begin{equation}
Y^2 = X^3 - g_2X - g_3 \label{Weierstrass}
\end{equation}
This is indeed possible and is accomplished in practice by performing a further birational transformation of the complex variable $\nu$ in the Legendre normal form (\ref{Legendre_n.f.}) of the equation
\begin{equation}
\nu = \frac{X-e_3}{e_1-e_3} \label{Legendre-Weierstrass}
\end{equation}
while also substituting $\mu=Y/\rho$. Under the transformation (\ref{Legendre-Weierstrass}) the abelian 1-form (\ref{omega_Legendre}) becomes
\begin{equation}
\varpi = \frac{dX}{2\sqrt{X^3 - g_2X - g_3}} \label{omega_Weierstrass}
\end{equation}
The discriminant of the Weierstrass cubic (\ref{Weierstrass})
\begin{equation}
\Delta = 4g_2^3 - 27 g_3^2 = \rho^6 k^4 k'^4 = [(e_1-e_2)(e_2-e_3)(e_3-e_1)]^2 
\end{equation}
is strictly positive as long as the elliptic modulus $k$ is not $0$ or $1$.

Let us introduce now in place of the Majorana coefficients the related real parameters
\begin{equation}
x_{\pm} = \frac{x\pm6|z|}{3} \qquad\qquad v_+ = \mbox{Im}\frac{v}{\sqrt{z}} \qquad\qquad v_- = \mbox{Re}\frac{v}{\sqrt{z}} \label{xpm_vpm}
\end{equation}
In terms of these, the invariants (\ref{g2_Majorana}) and (\ref{g3_Majorana}) can be written as follows
\begin{eqnarray}
g_2 &\hspace{-7pt} = \hspace{-7pt}& \ x_+^2+x_+x_-+x_-^2 + \frac{1}{4}(x_+\!-x_-)(v_-^2+v_+^2) \label{g2}   \\ 
g_3 &\hspace{-7pt} = \hspace{-7pt}& -(x_++x_-)x_+x_- - \frac{1}{4}(x_+\!-x_-)(x_+v_-^2+x_-v_+^2) \label{g3}
\end{eqnarray}
This form of the Weierstrass coefficients facilitates a key insight, namely that  the four points with $(X,Y)$-coordinates 
\begin{equation}
\begin{array}{lcr}
\displaystyle{(x_-,v_-(x_+\!-x_-)/2)} &\quad& \displaystyle{(x_+,iv_+(x_+\!-x_-)/2)} \\ [9pt]
\displaystyle{(x_-,v_-(x_-\!-x_+)/2)} &\quad& \displaystyle{(x_+,iv_+(x_-\!-x_+)/2)}
\end{array} \label{four_points}
\end{equation}
are points on the ${\cal O}(4)$ curve in the Weierstrass representation, {\it i.e.}, they satisfy the equation (\ref{Weierstrass}). This can be checked by direct substitution. Note that the pairs of points on the two colums in (\ref{four_points}) are conjugated under the elliptic involution. The pairs of points along the two diagonals are conjugated under the $\mathbb{Z}_2$ action
\begin{eqnarray}
x_-	& \longleftrightarrow & x_+ \nonumber \\
v_-	& \longleftrightarrow & iv_+ \label{Z2@O4}
\end{eqnarray}
Clearly, this action leaves the coefficients $g_2$ and $g_3$ invariant.

In the Legendre normal form of the curve the $\nu$-coordinates  corresponding to the points (\ref{four_points}) on the Weierstrass cubic are
\begin{equation}
\nu_{\pm} = \frac{x_{\pm} - e_3}{e_1-e_3} \label{nu_pm-x_pm}
\end{equation}
Using equations (\ref{Weierstrass_roots}), (\ref{k-modulus}) and the relations between the coefficients and the roots in (\ref{eta4}), we find that
\begin{eqnarray}
\nu_{\pm} 		&\hspace{-7pt}=\hspace{-7pt}& \frac{(1 \pm |\alpha\beta|)^2}{(1+|\alpha|^2)(1+|\beta|^2)} \nonumber \\
1-\nu_{\pm}	&\hspace{-7pt}=\hspace{-7pt}& \frac{(|\alpha| \mp |\beta|)^2}{(1+|\alpha|^2)(1+|\beta|^2)} \nonumber \\
k^2-\nu_{\pm}	&\hspace{-7pt}=\hspace{-7pt}& \frac{(\! \sqrt{\alpha\bar{\beta}} \mp \!\sqrt{\rule{0pt}{8pt}\bar{\alpha}\beta})^2}{(1+|\alpha|^2)(1+|\beta|^2)} \label{nu_pm}
\end{eqnarray}
A quick inspection of these relations yields the inequalities
\begin{equation}
0 < \nu_- < k^2 < \nu_+ < 1
\end{equation}
On the Weierstrass side, together with the obvious ordering of the Weierstrass roots, they imply that
\begin{equation}
e_3 < x_- < e_2 < x_+ < e_1 < e_0 = \infty \label{ineq_W}
\end{equation}

\subsection{Radial ${\cal O}(4)$ invariants} \label{O4-radii}

In the limit when an ${\cal O}(4)$ multiplet degenerates to the square of an ${\cal O}(2)$ multiplet, that is,  when $\eta^{(4)} \longrightarrow (\eta^{(2)})^2$, one can show that $\rho \longrightarrow |\vec{r}\,|^2$, where $\vec{r}$ is the $\mathbb{R}^3$ vector of the type (\ref{O2-radii}) associated to the ${\cal O}(2)$ multiplet. This suggests that, in the (generic, non-degenerate) ${\cal O}(4)$ case, positively-defined $SO(3)$-invariant quantities which are proportional to $\sqrt{\rho}$ could be interpreted as some sort of radii, too. There are essentially only two independent quantities satisfying these requirements, namely
\begin{equation}
r=\frac{1}{2\omega}>0 \qquad \mbox{and} \qquad r'=\frac{i\pi}{2\omega'}>0 \label{r&r'}
\end{equation}
The disparity in the normalization factors serves a subsequent formal purpose but is otherwise irrelevant.

These new invariants are not independent from the ones we have already introduced. One can  express  $g_2$ and $g_3$ in terms of $r$ and $r'$. The precise expressions follow for instance from the well-known double series representation formulas for the Weierstrass coefficients
\begin{eqnarray}
g_2 &\hspace{-7pt} = \hspace{-7pt}& 15\sum_{m,m'=-\infty}^{\infty} \hspace{-14pt}{}' \hspace{10pt} \frac{1}{(2m\omega+2m'\omega')^4}  \\ [2pt]
g_3 &\hspace{-7pt} = \hspace{-7pt}& 35\sum_{m,m'=-\infty}^{\infty} \hspace{-14pt}{}' \hspace{10pt} \frac{1}{(2m\omega+2m'\omega')^6}
\end{eqnarray}
where the prime sum symbol signifies that the term with $(m,m')=(0,0)$ must be omitted. Alternatively,  each of these double series can be recast as a single Lambert-type $q$-series
\begin{eqnarray}
g_2 &\hspace{-7pt} = \hspace{-7pt}& \hspace{3pt} \frac{1}{3}\, \left(\frac{\pi}{2\omega}\right)^4 \!\left( 1+240\sum_{n=1}^{\infty} n^3\frac{q^{2n}}{1-q^{2n}} \right) \\ [2pt]
g_3 &\hspace{-7pt} = \hspace{-7pt}& \frac{2}{27} \left(\frac{\pi}{2\omega}\right)^6 \!\left( 1-504\sum_{n=1}^{\infty} n^5\frac{q^{2n}}{1-q^{2n}} \right)
\end{eqnarray}
where $q=\exp(i\pi\tau)$ is the elliptic {\it nome} and $\tau=\omega'/\omega$ is the elliptic modulus. Since the Weierstrass coefficients are invariant under the modular transformation $\tau'=-1/\tau$, see for example the discussion in section \ref{mod-transf}, their $q'$-series expansions are formally identical, but with $q$ replaced by $q'$ and $\omega$ by $\omega'$.
In terms of the radii, 
\begin{equation}
q = e^{-\pi^2r/r'} \qquad\mbox{and}\qquad q' = e^{-r'/r} 
\end{equation} 
The fact that $r, r' > 0$ implies that $0<q,q'<1$, which in turn guarantees convergence. The two asymptotic regions $r>>r'$ and $r<<r'$ can be analyzed perturbatively by performing expansions in $q$ respectively $q'$.

\subsection{Euler-angle parametrization}

The rotational transformation properties of the Majorana coefficients suggest an alternative parametrization of the ${\cal O}(4)$ multiplet in which the rotational structure appears more explicitly. Consider the following {\it Ansatz}
\begin{eqnarray}
z &\hspace{-7pt} = \hspace{-7pt}& \sqrt{1} \sum_{m=-2}^2 \! D^{\, (2)}_{\!-2m}(\phi,\theta,\psi) \chi^2_m \nonumber \\
v &\hspace{-7pt} = \hspace{-7pt}& \sqrt{4} \sum_{m=-2}^2 \! D^{\, (2)}_{\!-1m}(\phi,\theta,\psi) \chi^2_m \nonumber \\
x &\hspace{-7pt} = \hspace{-7pt}& \sqrt{6} \sum_{m=-2}^2 \! D^{\, (2)}_{\!\phantom{-}0m}(\phi,\theta,\psi) \chi^2_m \label{spherical_param1}
\end{eqnarray}
Three of the new parameters will be the Euler angles $\phi$, $\theta$ and $\psi$. Wigner's rotation matrices ensure the right transformation properties. This leaves two rotation-invariant parameters on which the 5-component $\chi^2_m$ can depend. Observe now that due to the rotational invariance of $g_2$ and $g_3$ one can replace in equations (\ref{g2_Majorana}) and (\ref{g3_Majorana}) $z$, $v$ and $x$ by $\sqrt{1}\,\chi^2_{-2}$, $\sqrt{4}\,\chi^2_{-1}$ and $\sqrt{6}\,\chi^2_{0}$, respectively. Passing then to a form similar to (\ref{g2}) and (\ref{g3}), $x_{\pm}$ gets replaced by $e_{\pm} = (\sqrt{6}\chi^2_0 \pm 6 |\chi^2_{-2}|)/3$, and it is clear that if we put $\chi^2_{-1} = 0 = \chi^2_{+1}$ and compare the remainder with (\ref{G2}) and (\ref{G3}) we can identify $e_+$ and $e_-$ with any two of $e_1$, $e_2$ and $e_3$. Taking for instance $e_+=e_1$ and $e_-=e_3$, it follows that 
$\chi^2_0 = -\sqrt{6}\,e_2/4$ and $|\chi^2_{-2}| = (e_1-e_3)/4$. A choice consistent with these constraints would be, {\it e.g.},
\begin{equation}
\chi^2_{m=-2\cdots+2} = \frac{1}{4}
\left( \!\!
\begin{array}{c}
e_1-e_3 \\ 
0 \\ [2pt]
-\sqrt{6}\,e_2 \\ [2pt]
0 \\ 
e_1-e_3
\end{array}
\!\! \right) \label{spherical_param2}
\end{equation}
It depends, as required, on two rotation-independent parameters, since $e_1$, $e_2$ and $e_3$ satisfy (\ref{G1}). The parametrization of the ${\cal O}(4)$ multiplet given by (\ref{spherical_param1}) and (\ref{spherical_param2}) corresponds essentially to the one used by Ro\v{c}ek and Ivanov in \cite{Ivanov:1995cy}.

\subsection{${\cal O}(4)$ elliptic integrals} \label{integrals} \label{SEC:key-O4-int}

In the course of the generalized Legendre transform constructions involving ${\cal O}(4)$ multiplets that we shall encounter, the calculation of either the $F$-function contour integral or its first derivatives reduces  to the evaluation of the following set of integrals
\begin{equation}
{\cal I}_m = \int_{\Gamma} \frac{d\zeta}{\zeta} \frac{\zeta^m}{2\sqrt{\eta^{(4)}}} \label{I_m}
\end{equation}
with $\Gamma$ an integration contour which may be either open or closed, depending on context, and $m$ an integer taking values from $-2$ to $2$. In fact, it suffices to consider only $m = 0,1,2$, since the integrals corresponding to $m$ and $-m$ are complex conjugated to each other, modulo a shift. More precisely,
\begin{equation}
{\cal I}_{-m} = (-)^m \bar{\cal I}_m \pm 2\pi i \frac{m\bar{\beta}^{m-1}}{\sqrt{z}} \label{I_-m}
\end{equation}
This can be seen by changing in (\ref{I_m}) the integration variable $\zeta$ to $-1/\bar{\zeta}$ and then deforming the resulting contour back to the original one; in the process, one picks up a residue, which accounts for the shift term. Shifts will be discarded by means of a doubling trick: we can always choose two contours, one which gives a $+$ and one which gives a $-$ in (\ref{I_-m}); by summing the two contributions up, the residue terms will mutually cancel. 

In \cite{MM3} we evaluate these integrals in terms of Weierstrass elliptic functions by going to the Jacobian of the elliptic curve associated to the ${\cal O}(4)$ multiplet. We quote here the results for closed contour integrals which are of primary interest for us at this stage, and refer to \cite{MM3} for proofs. Integrating successively in the $\zeta$-plane along the closed contours $\Gamma_1$ surrounding the roots $-1/\bar{\beta}$ and $\alpha$, $\Gamma_2$ surrounding the roots $\alpha$ and $-1/\bar{\alpha}$ and $\Gamma_3$ surrounding the roots $-1/\bar{\alpha}$ and $-1/\bar{\beta}$, we obtain
\begin{eqnarray}
{\cal I}_0^{(i)} &\hspace{-7pt} = \hspace{-7pt}& 2\omega_i \label{I0-c} \\ [8pt]
{\cal I}_1^{(i)} &\hspace{-7pt} = \hspace{-7pt}& \frac{1}{\sqrt{z}}[\pi_i(x_+)+\pi_i(x_-)]  \label{I1-c} \\
{\cal I}_2^{(i)} &\hspace{-7pt} = \hspace{-7pt}& -\frac{1}{2z} \left[ 2\eta_i + (x_+\!+x_-)\omega_i - \frac{v}{\sqrt{z}}[\pi_i(x_+)+\pi_i(x_-)]  \right] \label{I2-c}
\end{eqnarray}
where the index $i$ takes values from $1$ to $3$, $\omega_1$, $\omega_2$, $\omega_3$ and $\eta_1$, $\eta_2$, $\eta_3$ are the Weierstrass half-periods and quasi-periods, respectively, and the $\pi$-functions are Jacobi elliptic integrals of the third kind; we give their precise definitions in the Appendix, where we also work out their properties.

\section{The Atiyah-Hitchin metric} \label{SEC:AH-metric}

\subsection{The hyperk\"{a}hler structure}

The $F$-function that generates the Atiyah-Hitchin metric in the generalized Legendre transform construction is given, according to \cite{Ivanov:1995cy}, by
\begin{equation}
F = F_2+F_1 =  \frac{1}{2\pi i h} \oint_{\Gamma_0} \frac{d\zeta}{\zeta} \eta^{(4)} - \oint_{\Gamma} \frac{d\zeta}{\zeta} \sqrt{\eta^{(4)}} \label{F_AH}
\end{equation}
$\Gamma_0$ is an integration contour around $\zeta=0$ whereas $\Gamma$ is a contour that winds once around the branch-cut between the roots $\alpha$ and $-1/\bar{\beta}$; $h$  is a constant coupling scale.

The first integral can be evaluated by means of a straightforward application of Cauchy's integral formula. One gets
\begin{equation}
F_2 = \frac{x}{h} \label{F2}
\end{equation}
The evaluation of the second integral presents a more challenging problem. Observe that it satisfies the homogeneity property
\begin{equation}
F_1 = 2\bar{z}\frac{\partial F_1}{\partial\bar{z}} + 2\bar{v}\frac{\partial F_1}{\partial\bar{v}} + 2x\frac{\partial F_1}{\partial x} + 2v\frac{\partial F_1}{\partial v} + 2z\frac{\partial F_1}{\partial z} \label{F1-hom}
\end{equation}
and that each derivative is equal to a complete integral of the type (\ref{I_m}), which we have already evaluated. More precisely, 
\begin{equation}
\frac{\partial F_1}{\partial z} = - {\cal I}_2^{(1)} 
\qquad\qquad 
\frac{\partial F_1}{\partial v} = {\cal I}_1^{(1)} 
\qquad\qquad 
\frac{\partial F_1}{\partial x} = - {\cal I}_0^{(1)} \label{dF=I_m}
\end{equation}
with ${\cal I}_{m = 0,1,2}^{(1)}$ given in (\ref{I0-c}) - (\ref{I2-c}).  To obtain a real $F_1$ we should consider a combination of contours around the branch-cuts between the roots $\alpha$ and $-1/\bar{\beta}$ as well as between $\beta$ and $-1/\bar{\alpha}$ in such a way as to cancel the residue terms in (\ref{I_-m}). Henceforth, we will automatically assume that this is done, and we will just ignore them. With this assumption, equations (\ref{F1-hom}), (\ref{dF=I_m}) and (\ref{I0-c}) - (\ref{I2-c}) yield
\begin{equation}
F_1 	= 4\eta - 4(x_+\!+x_-)\omega + 2 v_-\pi(x_-) + 2iv_+\pi(x_+) \label{F1}
\end{equation}
From the equations (\ref{F2}) and (\ref{dF=I_m}) one also readily obtains the derivatives
\begin{eqnarray}
\frac{\partial F}{\partial v} &\hspace{-7pt} = \hspace{-7pt}& \frac{\pi(x_-) + \pi(x_+)}{\sqrt{z}} \label{Fv} \\
\frac{\partial F}{\partial x} &\hspace{-7pt} = \hspace{-7pt}& \frac{1}{h} - 2\omega \label{Fx}
\end{eqnarray}
We are now prepared to impose the generalized Legendre relations, which in this case read
\begin{eqnarray}
\frac{\partial F}{\partial v} &\hspace{-7pt} = \hspace{-7pt}& u \label{Lrel1} \\ [2pt]
\frac{\partial F}{\partial x} &\hspace{-7pt} = \hspace{-7pt}& 0 \label{Lrel2}
\end{eqnarray}
Equation (\ref{Lrel2}) defines implicitly the Atiyah-Hitchin manifold as a codimension 1 subspace in the five real-dimensional space of moduli of $\eta^{(4)}$ sections. It takes the remarkably simple form 
\begin{equation}
2\omega = 1/h \label{AH}
\end{equation}
Equation (\ref{Lrel1}), on the other hand, serves to introduce a second holomorphic coordinate $u$ (the other one being $z$). Note that due to the inequalities (\ref{ineq_W}), $\pi(x_-)$ is real whereas $\pi(x_+)$ is purely imaginary. Then equations (\ref{Fv}) and (\ref{Lrel1}) imply that $\pi(x_-) = \mbox{Re}\, u\sqrt{z}$ and $ \pi(x_+) = i\, \mbox{Im}\, u\sqrt{z}$, and so
\begin{equation}
2 v_-\pi(x_-) + 2iv_+\pi(x_+) = uv + \overline{uv} \label{Leg_tr}
\end{equation}
{\it It follows that the Jacobi terms in $F$, in conjunction with the Legendre relation that introduces the second holomorphic coordinate, play a pivotal role in canceling the quadratic terms in the Legendre transform  that yields the K\"{a}hler potential.}  This constitutes a generic mechanism for ${\cal O}(4)$ multiplets. Putting together equations~(\ref{F1}), (\ref{F2}) and (\ref{Leg_tr}) and making also use of equation~(\ref{AH}), one obtains quite effortlessly a formula for the K\"{a}hler potential corresponding to the complex structure to which the holomorphic coordinates $z$ and $u$ are associated, {\it i.e.},
\begin{equation}
K = 4\eta -(x_+\!+x_-)\omega \label{Kahler_AH}
\end{equation}
A K\"{a}hler potential for the Atiyah-Hitchin metric has also been derived by Olivier \cite{Olivier:1991pa}, who followed a different, symmetry-based approach. It is straightforward to check that the rotation-invariant term of (\ref{Kahler_AH}), namely $4\eta$, coincides, up to a constant factor, with the corresponding part of Olivier's potential.

We found {\it a posteriori} that the metric takes a simpler form if one uses instead of $u$ and $z$ a new pair of holomorphic variables defined as follows
\begin{equation}
U = u\sqrt{z} \qquad\quad Z = 2\sqrt{z} \label{U_and_Z}
\end{equation}
This transformation is a holomorphic symplectomorphism, it leaves the hyperk\"{a}hler holomorphic (2,0)-form invariant, 
\begin{equation}
\omega^{+} = dZ \wedge dU
\end{equation}
Taking also into account the reality properties of $\pi(x_{\pm})$, equation (\ref{Lrel1}) is then equivalent to
\begin{equation}
U = \pi(x_-)+\pi(x_+) \qquad\quad \bar{U} = \pi(x_-)-\pi(x_+) \label{U}
\end{equation} 

To derive the metric one can follow two equivalent paths:  one is to use the general formula (\ref{GLT-metric}) giving the metric in terms of the second derivatives of $F$ without resorting explicitly to holomorphic variables, the other is to use the K\"{a}hler potential (\ref{Kahler_AH}) and the manifest holomorphic structure. 
In the first approach, the second derivatives of $F$ follow from derivating once more $F_v$ and $F_x$, given by equations~(\ref{Fv}) and (\ref{Fx}). The absence of $F_z$ from this list can be compensated by the relations (\ref{diff-eqs}) between the second derivatives. Regarding $\omega$ and $\pi(x_{\pm})$ as functions of $g_2$, $g_3$ and (in the latter case only) $x_{\pm}$, one makes use of the differentiation formulas (\ref{d_omega}) and (\ref{d_piX}) from the Appendix in combination with the explicit forms (\ref{g2_Majorana}) and (\ref{g3_Majorana}) of $g_2$ and $g_3$.
In the second approach, one employs equations~(\ref{AH}) and (\ref{U}), regarding now $\pi(x_{\pm})$ as a function of $\omega$, $\eta$ and $x_{\pm}$. From the equations~(\ref{U}) together with $|Z|^2 = x_+\!-x_-$, by making use of the differentiation relation (\ref{d_piX2}) one can obtain the partial derivatives of $x_{\pm}$ and $\eta$ with respect to the holomorphic variables $U$ and $Z$,
\begin{eqnarray}
d\eta &\hspace{-7pt} = \hspace{-7pt}& \frac{(A_+\!-A_-)dU+(A_+\!+A_-)d\bar{U}+2A_+A_-(\bar{Z}dZ+Zd\bar{Z})}{2(A_-B_+\!-A_+B_-)} \\
dx_{\pm} &\hspace{-7pt} = \hspace{-7pt}& \frac{(B_+\!-B_-)dU+(B_+\!+B_-)d\bar{U}+2A_{\mp}B_{\pm}(\bar{Z}dZ+Zd\bar{Z})}{2(A_-B_+\!-A_+B_-)}
\end{eqnarray}
where, with $y_+ = iv_+(x_+\!-x_-)/2$ and $y_- = v_-(x_-\!-x_+)/2$ as in (\ref{four_points}) and $V$ given in (\ref{vee}),
\begin{eqnarray}
A_{\pm} &\hspace{-7pt} = \hspace{-7pt}& \frac{\eta+x_{\pm}\omega}{2y_{\pm}} \nonumber \\
B_{\pm} &\hspace{-7pt} = \hspace{-7pt}& \frac{V\omega+x_{\pm}}{y_{\pm}}
\end{eqnarray}
The components of the metric can now be computed directly from the K\"{a}hler potential (\ref{Kahler_AH}). Following either path, one obtains
\begin{eqnarray}
K_{Z\bar{Z}} &\hspace{-7pt} = \hspace{-7pt}& \phantom{+}\frac{1}{{\cal Q}|Z|^2}{\cal K}_4  \nonumber \\
K_{U\bar{Z}} &\hspace{-7pt} = \hspace{-7pt}& -\frac{1}{2{\cal Q}\bar{Z}} ( v_-{\cal K}_{3+}  +  i v_+{\cal K}_{3-} )  \nonumber \\
K_{Z\bar{U}} &\hspace{-7pt} = \hspace{-7pt}& -\frac{1}{2{\cal Q}Z} (  v_-{\cal K}_{3+} - iv_+{\cal K}_{3-} ) \nonumber \\
K_{U\bar{U}} &\hspace{-7pt} = \hspace{-7pt}& \phantom{+}\frac{1}{{\cal Q}|Z|^2}{\cal K}_2 \label{AH_metric_comps}
\end{eqnarray}
where
\begin{equation}
{\cal Q} = (\eta+e_1\omega) (\eta+e_2\omega)(\eta+e_3\omega) 
\end{equation}
and
\begin{eqnarray}
{\cal K}_2 		&\hspace{-7pt} = \hspace{-7pt}& (g_2-3 x_+x_-) \eta^2 - [6g_3+2(x_+\!+x_-)g_2] \omega\eta  + [g_2^2+3(x_+\!+x_-)g_3+x_+x_-g_2]\omega^2 \nonumber \\ [4pt]
{\cal K}_{3} 	&\hspace{-7pt} = \hspace{-7pt}& \eta^3 +3x_{\pm} \omega\eta^2 + g_2\, \omega^2\eta - (2g_3+x_{\pm}g_2)\omega^3 \\ [4pt]
 {\cal K}_4 	&\hspace{-7pt} = \hspace{-7pt}& \eta^4 + 2(x_+\!+x_-) \omega\eta^3  + (g_2+3x_+x_-) \omega^2\eta^2 - 2g_3 \omega^3\eta - [(x_+\!+x_-)g_3+x_+x_-g_2]\omega^4 \nonumber
\end{eqnarray}
One verifies that the following Monge-Amp\`{e}re equation holds 
\begin{equation}
\det K_{(Z,U)} = 1
\end{equation}
Putting things together, the metric eventually reads
\begin{equation}
ds^2 = \phi\, dU d\bar{U} + \phi^{-1}\! \left(dZ + {\cal A}\right)\left(d\bar{Z} + \bar{{\cal A}}\right) \label{AH_metric_hol}
\end{equation}
with
\begin{equation}
\phi = \frac{\cal Q}{{\cal K}_4}|Z|^2 \qquad\mbox{and}\qquad {\cal A} = -\frac{v_-{\cal K}_{3+} \!+ iv_+{\cal K}_{3-}}{2{\cal K}_4} ZdU
\end{equation}
It has three isometries induced by the natural $SO(3)$ action on the parameter space of the ${\cal O}(4)$ section, generated by the three vector fields
\begin{eqnarray}
 X_{3\, } &\hspace{-7pt} = \hspace{-7pt}&  -2i \!\left(\!Z\frac{\partial}{\partial Z} - \bar{Z}\frac{\partial}{\partial \bar{Z}}\right) \\
X_- &\hspace{-7pt} = \hspace{-7pt}& (v_-\!+iv_+)\frac{\partial}{\partial Z}  + \frac{(\eta+x_-\omega)+(\eta+x_+\omega)}{Z}\frac{\partial}{\partial U} + \frac{(\eta+x_-\omega)-(\eta+x_+\omega)}{Z}\frac{\partial}{\partial \bar{U}} \hspace{20pt} \\ [6pt]
X_+ &\hspace{-7pt} = \hspace{-7pt}& \overline{X_-}
\end{eqnarray}
These preserve the metric, form an $SO(3)$ algebra,\footnote{We use the standard conventions $X_{\pm} = (X_1 \pm i X_2)/2$ and $\omega^{\pm} = (\omega_1\pm i\omega_2)/2$.}
\begin{equation}
[X_i,X_j] = 2\epsilon_{ijk\,} X_k \label{isom_SO(3)_algebra}
\end{equation}
and rotate the hyperk\"{a}hler structure
\begin{equation}
{\cal L}_{X_i} \omega_j = 2\epsilon_{ijk\,} \omega_k \label{K_str_rot}
\end{equation}
From the last property it is clear that the actions of the vector fields $X_i$ are hamiltonian with respect to the corresponding K\"{a}hler 2-forms $\omega_i$ and thus, provided that the manifold is simply connected, they must have associated moment maps. We get, explicitly,
\begin{eqnarray}
\mu_1 &\hspace{-7pt} = \hspace{-7pt}& 4\eta-(x_+\!+x_-)\omega + {\rm Re} Z^2 \omega \label{mu1} \\ [4pt]
\mu_2 &\hspace{-7pt} = \hspace{-7pt}& 4\eta-(x_+\!+x_-)\omega - {\rm Re} Z^2 \omega \label{mu2} \\ [4pt]
\mu_3 &\hspace{-7pt} = \hspace{-7pt}& 4\eta+2(x_+\!+x_-)\omega \label{mu3}
\end{eqnarray}
In \cite{Hitchin:1986ea,MR1726926} it is shown that if on a simply-connected hyperk\"{a}hler manifold there exists a vector field $X_i$ which acts on the hyperk\"{a}hler structure as in (\ref{K_str_rot}) then the corresponding moment map is a K\"{a}hler potential for any complex structure orthogonal\footnote{A hyperk\"{a}hler manifold with standard complex structures $I$, $J$, $K$ has in fact a 2-sphere worth of complex structures compatible with the metric and the Levi-Civita connection, given by $xI+yJ+zK$ for all unit vectors $(x,y,z)$  from $\mathbb{R}^3$. One refers to two complex structures as being orthogonal when the corresponding unit vectors are orthogonal with respect to the scalar product on $\mathbb{R}^3$.} to the one preserved by $X_i$. In agreement with this theorem, the moment maps (\ref{mu1}) and (\ref{mu2}) belong manifestly to the K\"{a}hler class of (\ref{Kahler_AH}).

In the derivation of the relations (\ref{isom_SO(3)_algebra}) through (\ref{mu3}) we found the following partial results useful
\begin{eqnarray}
Z\frac{\partial K}{\partial Z}  &\hspace{-7pt} = \hspace{-7pt}& \bar{Z}\frac{\partial K}{\partial\bar{Z}} =  \frac{\mu_3}{2} \\ 
Z\frac{\partial\mu_3}{\partial Z}  &\hspace{-7pt} = \hspace{-7pt}& \bar{Z}\frac{\partial\mu_3}{\partial\bar{Z}}  = \frac{2{\cal K}_4}{\cal Q} \\
\frac{\partial K}{\partial U} &\hspace{-7pt} = \hspace{-7pt}& -(v_-\!+iv_+) \\ 
\frac{\partial\mu_3}{\partial U} &\hspace{-7pt} = \hspace{-7pt}& -\frac{v_-{\cal K}_{3+} \!+ iv_+{\cal K}_{3-}}{{\cal Q}}
\end{eqnarray}

\subsection{The $SO(3)$ structure}

The connection between the complex holomorphic basis form (\ref{AH_metric_hol}) of the metric and the well-known form of Atiyah and Hitchin \cite{Atiyah:1988jp}, in which the (hyper-)complex structure is obscure but the non-triholomorphic $SO(3)$ isometry is manifest, emerges through a change of variables. Specifically, one needs to switch from the complex holomorphic coordinate basis $Z$, $\bar{Z}$, $U$, $\bar{U}$ to the coordinate basis given by the $SO(3)$ angles $\phi$, $\theta$, $\psi$ and the elliptic {\it nome} $q$. The Majorana coefficients $z$, $\bar{z}$, $v$, $\bar{v}$ and $x$ will serve as intermediate variables. The Jacobian matrix of this transformation is computed as follows:  the partial derivatives of $Z$ and $\bar{Z}$ with respect to the Majorana coefficients are quite trivially computed from the second relation (\ref{U_and_Z}); the partial derivatives of $U$ and $\bar{U}$ on the other hand can be obtained from (\ref{U}), regarding $\pi(x_{\pm})$ as functions of $g_2$, $g_3$ and $x_{\pm}$. In this way one gets, {\it e.g.},
\begin{equation}
\frac{\partial U}{\partial z} = \frac{\bar{Z}}{4Z}(M \eta^* + N\omega^*) \label{dUdz}
\end{equation}
with $\eta^*$ and $\omega^*$ defined in (\ref{omet1}) and
\begin{eqnarray}
M &\hspace{-7pt} = \hspace{-7pt}& v_-(5x_+\!-x_-)-iv_+(5x_-\!-x_+)  \nonumber \\ [4pt]
N &\hspace{-7pt} = \hspace{-7pt}& v_-(8g_2-3x_+x_-\!-9x_+^2) - iv_+(8g_2-3x_+x_-\!-9x_-^2) \label{M_and_N}
\end{eqnarray}
For the remaining derivatives, a detour through the relations (\ref{diff-eqs}) allows for simpler calculations. Specifically, 
\begin{eqnarray}
\frac{\partial U}{\partial x} &\hspace{-7pt} = \hspace{-7pt}& \frac{Z}{2} F_{xv} = \phantom{+}\frac{Z}{2} F_{vx} =  - Z \frac{\partial \omega}{\partial v} \nonumber \\
\frac{\partial U}{\partial \bar{v}} &\hspace{-7pt} = \hspace{-7pt}& \frac{Z}{2} F_{\bar{v}v} = -\frac{Z}{2} F_{xx} = \phantom{+} Z \frac{\partial \omega}{\partial x} \nonumber \\
\frac{\partial U}{\partial v} &\hspace{-7pt} = \hspace{-7pt}& \frac{Z}{2} F_{vv} = \phantom{+}\frac{Z}{2} F_{zx} = - Z \frac{\partial \omega}{\partial z} \nonumber \\
\frac{\partial U}{\partial \bar{z}} &\hspace{-7pt} = \hspace{-7pt}& \frac{Z}{2} F_{\bar{z}v} = - \frac{Z}{2} F_{\bar{v}x} = \phantom{+} Z \frac{\partial \omega}{\partial \bar{v}} \label{dw_1} 
\end{eqnarray}
\noindent Explicit forms along the lines of (\ref{dUdz})-(\ref{M_and_N}) can be easily obtained by means of the formula (\ref{d_omega}) in combination with the equations (\ref{g2_Majorana}) and (\ref{g3_Majorana}). It remains now to compute the derivatives of the Majorana coefficients with respect to $\phi$, $\theta$, $\psi$ and $q$. To do this, we use the equations (\ref{spherical_param1}) and (\ref{spherical_param2}) together with (\ref{W_roots_theta}) to express the Majorana coefficients in terms of $\phi$, $\theta$, $\psi$, $q$ and $\omega$. The half-period $\omega$ is fixed by the equation (\ref{AH}). The derivatives may then be computed with the help of, among other things, the formulas (\ref{dW_roots_dtau}), which gives the partial derivatives of the Weierstrass roots with respect to the {\it nome}. 

In the new coordinate basis the metric takes the form\footnote{An overall numerical factor is omitted.}
\begin{equation}
ds^2 \sim (abc)^2  \frac{dq^2}{q^2} + a^2\sigma_1^2 + b^2\sigma_2^2 + c^2\sigma_3^2
\end{equation}
with $\sigma_1$, $\sigma_2$ and $\sigma_3$ the (left-)invariant Cartan-Maurer 1-forms of $SO(3)$ and $a$, $b$, $c$ determined by
\begin{eqnarray}
2bc &\hspace{-7pt} = \hspace{-7pt}& -\left(\frac{2}{\pi}\right)^2\! \omega (\eta+e_1 \omega) \nonumber \\ [0pt]
2ab &\hspace{-7pt} = \hspace{-7pt}& -\left(\frac{2}{\pi}\right)^2\!\omega (\eta+e_2 \omega) \nonumber \\ [0pt]
2ca &\hspace{-7pt} = \hspace{-7pt}& -\left(\frac{2}{\pi}\right)^2\!\omega (\eta+e_3 \omega) 
\label{bc_ab_ca}
\end{eqnarray}
Comparison with the theta-function representation formula (\ref{lola}) from the Appendix yields eventually for the metric coefficients the expressions
\begin{eqnarray}
2bc &\hspace{-7pt} = \hspace{-7pt}&  \frac{\vartheta''_2(0,q)}{\vartheta_2(0,q)} \nonumber \\ 
2ab &\hspace{-7pt} = \hspace{-7pt}&  \frac{\vartheta''_3(0,q)}{\vartheta_3(0,q)} \nonumber \\ 
2ca &\hspace{-7pt} = \hspace{-7pt}&  \frac{\vartheta''_4(0,q)}{\vartheta_4(0,q)} 
\label{bc_ab_ca-theta}
\end{eqnarray}
A change of variable from the {\it nome} $q$ to the complementary {\it nome} $q'$ through the modular transformation (\ref{tau-prime}) gives the following alternative expression for the metric
\begin{equation}
ds^2 \sim \left(\frac{\ln q'}{\pi}\right)^2\! \left[(ABC)^2  \frac{dq^{\prime 2}}{q^{\prime 2}} + A^2\sigma_1^2 + B^2\sigma_2^2 + C^2\sigma_3^2\right]
\end{equation}
with the coefficients $A$, $B$, $C$ determined this time by
\begin{eqnarray}
2BC &\hspace{-7pt} = \hspace{-7pt}& \frac{2}{\ln q'} - \frac{\vartheta''_4(0,q')}{\vartheta_4(0,q')} \nonumber \\ 
2AB &\hspace{-7pt} = \hspace{-7pt}& \frac{2}{\ln q'} - \frac{\vartheta''_3(0,q')}{\vartheta_3(0,q')}  \nonumber \\ 
2CA &\hspace{-7pt} = \hspace{-7pt}& \frac{2}{\ln q'} - \frac{\vartheta''_2(0,q')}{\vartheta_2(0,q')} 
\label{BC_AB_CA-theta}
\end{eqnarray}
One can check that these two expressions correspond precisely to the metric of Atiyah and Hitchin, up to an irrelevant numerical scale factor.

\subsection{Large and small monopole separation limits}

In terms of the radial invariants introduced in (\ref{r&r'}), the generalized Legendre transform equation (\ref{AH}) reads
\begin{equation}
r = h = \mbox{const.}
\end{equation}
The obvious interpretation is that a constant distance scale is thus set into the problem. The remaining radius, $r'$,  takes the meaning of monopole separation distance when one views the Atiyah-Hitchin manifold as the moduli space of centered two-monopoles. As discussed in section \ref{O4-radii}, the elliptic {\it nome} and complementary {\it nome} are given by
\begin{equation}
q = e^{-\pi^2h/r'} \qquad\mbox{and}\qquad q' = e^{-r'/h}
\end{equation}
The small monopole separation limit $r'<<h$ thus corresponds to $q \rightarrow 0$, the large monopole separation limit $r'>>h$ to $q' \rightarrow 0$. These asymptotic regions can therefore be probed by series-expanding in $q$ respectively $q'$ the theta-function expressions (\ref{bc_ab_ca-theta}) respectively (\ref{BC_AB_CA-theta}) of the coefficients of the metric. The infinite series representation formulas (\ref{q-series}) make this a straighforward task.

A $q'$-series expansion gives
\begin{eqnarray}
a^2 &\hspace{-7pt} = \hspace{-7pt}& \frac{r'(r'-2h)}{2h^2} - \frac{4r^{\prime 2}(r'-2h)}{h^3}e^{-r'/h} + \frac{4r^{\prime 2}(4r^{\prime 2}-8hr'+h^2)}{h^4}e^{-2r'/h} + \cdots \nonumber \\
b^2 &\hspace{-7pt} = \hspace{-7pt}& \frac{2r'}{r'-2h} - \frac{16r^{\prime 2}(2r^{\prime 2}-6hr'+5h^2)}{h^2(r'-2h)^2}e^{-2r'/h} + \cdots \nonumber \\ 
c^2 &\hspace{-7pt} = \hspace{-7pt}& \frac{r'(r'-2h)}{2h^2} + \frac{4r^{\prime 2}(r'-2h)}{h^3}e^{-r'/h} + \frac{4r^{\prime 2}(4r^{\prime 2}-8hr'+h^2)}{h^4}e^{-2r'/h} + \cdots
\end{eqnarray}
Retaining only the non-exponential terms yields the asymptotic form of the Atiyah-Hitchin metric
\begin{equation}
ds_{\infty}^2 = h^{-1} \!\left[ \left(\frac{1}{2h} - \frac{1}{r'}\right) (dr^{\prime 2} + r^{\prime 2}\sigma_1^2 + r^{\prime 2}\sigma_3^2) +  \left(\frac{1}{2h} - \frac{1}{r'}\right)^{\!-1} \!\sigma_2^2 \right]
\end{equation}
This is a Euclidean Taub-NUT metric with negative mass parameter. It has a singularity at $r' = 2h$, far away from the asymptotic region, and thus harmless. 

On the other hand, expanding in the {\it nome} $q$, we get
\begin{eqnarray}
a^2 &\hspace{-7pt} = \hspace{-7pt}& 32\pi^2(e^{-2\pi^2h/r'}-4e^{-4\pi^2h/r'}+\cdots) \nonumber \\ [6pt]
b^2 &\hspace{-7pt} = \hspace{-7pt}& \frac{\pi^2}{2\ }(1-4e^{-\pi^2h/r'}+16e^{-2\pi^2h/r'} + \cdots) \nonumber \\
c^2 &\hspace{-7pt} = \hspace{-7pt}& \frac{\pi^2}{2\ }(1+4e^{-\pi^2h/r'}+16e^{-2\pi^2h/r'} + \cdots)
\end{eqnarray}
Truncating to order $q^2$ and changing the radial variable to $R = 4\, e^{-\pi^2h/r'}$, we obtain the small monopole separation limit of the metric
\begin{equation}
ds_0^2 = \frac{\pi^2}{2\ } [dR^2 + 4R^2\sigma_1^2 + (1-R+R^2)\sigma_2^2 + (1+R+R^2)\sigma_3^2]
\end{equation}

These limits of the Atiyah-Hitchin metric and the first few exponential corrections to them have been studied in \cite{Atiyah:1988jp,Gibbons:1995yw}, to which we refer for more details. Note that, unlike in \cite{Gibbons:1995yw}, here we do not need to choose the gauge  $f = -b/r'$  instead of the more symmetric $f = abc$ for the radial diagonal component of the metric. A set of expressions closely related to the form (\ref{bc_ab_ca-theta}) of the metric coefficients has been obtained in \cite{Hanany:2000fw} through solving a Halphen system of differential equations, see also \cite{Takhtajan:1992qb}. The equations (\ref{bc_ab_ca-theta}) and (\ref{BC_AB_CA-theta}) together with the series expansions (\ref{q-series}) allow for a straightforward computation of the corrections to both the large and the small separation limit of the metric virtually to any order.




\section{APPENDIX: Elliptic functions and integrals}

\subsection{Legendre and Weierstrass elliptic integrals}

The incomplete elliptic integrals of first, second and third kind in Legendre normal form are
\begin{eqnarray}
F(z,k) &\hspace{-7pt} = \hspace{-7pt}& \int_0^z \frac{dt}{\sqrt{(1-t^2)(1-k^2t^2)}} \\ [4pt]
E(z,k) &\hspace{-7pt} = \hspace{-7pt}& \int_0^z \sqrt{\frac{1-k^2t^2}{1-t^2}} dt \label{E_Leg} \\ [4pt]
\Pi(z,\nu,k) &\hspace{-7pt} = \hspace{-7pt}& \int_0^z\frac{1}{1-\nu t^2} \frac{dt}{\sqrt{(1-t^2)(1-k^2t^2)}}
\end{eqnarray}
In place of Legendre's integral of third kind, Jacobi introduced a modified version, namely
\begin{equation}
\Pi_{\cal J}(z,\nu,k) = \sqrt{\frac{(\nu-1)(\nu-k^2)}{\nu}} [\Pi(z,\nu,k)-F(z,k)]  \label{Pi_Jac}
\end{equation}
Jacobi's integral enjoys a great many formal advantages over Legendre's, for a detailed discussion in the Legendre formalism frame see {\it e.g.} \cite{MR0124532}.
The parameters $z$, $k$ and $\nu$ are termed {\it amplitude}, {\it modulus} and {\it characteristic}, respectively. The corresponding complete integrals are obtained by putting the amplitude equal to $1$: $K(k) = F(1,k)$, $E(k) = E(1,k)$, $\Pi(\nu,k) = \Pi(1,\nu,k)$ and $\Pi_{\cal J}(\nu,K(k)) = \Pi_{\cal J}(1,\nu,k)$ (we will justify later on this last notation).  $K(k)$ arises from integrating the elliptic abelian  differential form around one of the canonical cycles of the Legendre elliptic curve. The integral around the other canonical cycle yields $iK(k')$, with the {\it complementary modulus} $k'$ satisfying $k^2+k'^2=1$. The simplified notations $K(k)=K$, $K(k')=K'$, {\it etc.} are sometimes used.

In the Weierstrass theory the role of the incomplete elliptic integrals is played by 
\begin{eqnarray}
&& \phantom{+} \int \frac{dX}{2Y} = u + {\cal C} \label{Wei1} \\ [3pt]
&& -\int X\frac{dX}{2Y} = \zeta(u) + {\cal C} \label{Wei2} \\ [2pt]
&& -\int \frac{Y_0}{X-X_0} \frac{dX}{2Y} = \frac{1}{2} \ln \frac{\sigma(u+u_0)}{\sigma(u-u_0)} - u\,\zeta(u_0) + {\cal C} \label{pi_W} \label{Wei3}
\end{eqnarray}
where ${\cal C}$ is an indefinite integration constant, $(X,Y)$ and  $(X_0,Y_0)$ are points on the Weierstrass curve $Y^2 = X^3-g_2X-g_3$, $u$ and $u_0$ are the corresponding points on the Jacobian variety, and $\sigma(u)$, $\zeta(u)$ are the Weierstrass sigma respectively zeta pseudo-elliptic functions. The expressions on the r.h.s. are obtained by substituting $X$ and $Y$ with the corresponding Weierstrass elliptic functions, {\it i.e.},
\begin{equation}
X = \wp(u;4g_2,4g_3) \hspace{30pt}  2Y = \wp'(u;4g_2,4g_3) \label{inv-AJ-map}
\end{equation}
The derivation of the first two expressions is fairly straightforward and standard. The derivation of the third one requires the use of a variant of the addition theorem of the Weierstrass $\zeta$-function. 

The corresponding complete integrals are obtained by integrating in the complex $X$-plane along the closed countours $\Gamma_1$, surrounding the roots $e_2$ and $e_3$, $\Gamma_2$, surrounding the roots $e_3$ and $e_2$ and $\Gamma_3$, surrounding the roots $e_2$ and $e_1$, or, more precisely, on the Jacobian, from $u = \omega_2$ to $-\omega_3$, from $u = \omega_3$ to $-\omega_2$ and from $u = \omega_2$ to $-\omega_1$, respectively. We get
\begin{eqnarray}
&& \phantom{+} \oint_{\Gamma_i} \frac{dX}{2Y} = 2\omega_i \label{Wcei1} \\ [2pt]
&& -\oint_{\Gamma_i} X\frac{dX}{2Y} = 2\eta_i \label{Wcei2} \\ [2pt]
&& -\oint_{\Gamma_i} \frac{Y_0}{X-X_0} \frac{dX}{2Y}  = 2\
\begin{array}{|cc|}
u_0		& \omega_i \\
\zeta(u_0)	\!&\! \zeta(\omega_i) 
\end{array} 
\stackrel{\rm def}{=} 2\,\pi_i(X_0) \label{Wcei3}
\end{eqnarray}
where $u_0$ is the image of $(X_0,Y_0)$ through the Abel-Jacobi map and $i = 1,2,3$. Equation (\ref{Wcei3}) follows by way of the monodromy property of the Weierstrass $\sigma$-function in the r.h.s. of (\ref{pi_W}). The notation $\pi_i(X_0)$ is not quite rigorous, a more appropriate one would be for instance $\pi_i(X_0,Y_0)$ or $\pi_i(u_0)$. Nonetheless, for simplicity reasons as well as for other practical reasons soon to become clear, we use it in this form, but with the implicit {\it caveat} that it conceals a sign ambiguity. Clearly, only two out of three integrals of each set of integrals are independent, as $\omega_1+\omega_2+\omega_3 = 0$, $\eta_1+\eta_2+\eta_3 = 0$ and $\pi_1(X)+\pi_2(X)+\pi_3(X) = 0$. In line with the usual notation conventions $\omega_1 = \omega$, $\omega_3 = \omega'$, $\eta_1 = \eta$, $\eta_3 = \eta'$ we also denote $\pi_1(X) = \pi(X)$ and $\pi_3(X) = \pi'(X)$. The integrals $\omega$ and $\omega'$ respectively $\eta$ and $\eta'$  are termed half-periods and half-quasi-periods because, for any $m, m' \in \mathbb{Z}$
\begin{eqnarray}
\wp(u+2m\omega+2m'\omega') &\hspace{-7pt} = \hspace{-7pt}& \wp(u) \\ [4pt]
\zeta(u+2m\omega+2m'\omega') &\hspace{-7pt} = \hspace{-7pt}& \zeta(u) + 2m\eta+2m'\eta'
\end{eqnarray}
Based on Legendre's identity
\begin{equation}
\begin{array}{|cc|}
\omega'		& \omega \\
\zeta(\omega')	& \zeta(\omega) 
\end{array} 
=i\frac{\pi}{2} \label{Legendre-id}
\end{equation}
one determines that
\begin{eqnarray}
u \longrightarrow u + 2m\omega+2m'\omega' \qquad \Longrightarrow \qquad \pi(X) &\hspace{-7pt}\longrightarrow \hspace{-7pt}& \pi(X)\, + i\pi m'  \\ [4pt]
 \pi'(X) &\hspace{-7pt}\longrightarrow \hspace{-7pt}& \pi'(X) - i\pi m \nonumber
\end{eqnarray}
This means that $\pi(X)$ and $\pi'(X)$ are {\it not} elliptic functions. Legendre's identity can be also used to show that
\begin{equation}
\begin{array}{|cc|}
\omega'		& \omega \\
\pi'(X)	 	\!&\! \pi(X) 
\end{array} 
=i\frac{\pi}{2} u \label{Legendre-pi}
\end{equation}
for any $X=\wp(u)$.

The connection between the complete integrals in the Weierstrass and the Jacobi theories is given by the following formulas
\begin{eqnarray}
K &\hspace{-7pt} = \hspace{-7pt}&\sqrt{\rho}\, \omega \label{K_JW} \\ [7pt]
E &\hspace{-7pt} = \hspace{-7pt}&\frac{\eta+e_1\omega}{\sqrt{\rho}} \label{E_JW} \\
\Pi_{\cal J}(\nu,K) &\hspace{-7pt} = \hspace{-7pt}&\pi(X)+i\frac{\pi}{2} \label{Pi_JW}
\end{eqnarray}
where $X$ is related to $\nu$ as in equation (\ref{Legendre-Weierstrass}). Similar relations hold for the corresponding primed quantities. Equations (\ref{E_JW}) and (\ref{Pi_JW}) can be proved by performing the following changes of integration variable in (\ref{E_Leg}) and (\ref{Pi_Jac}): $t^2=\tilde{\nu}/k^2$ respectively $t^2=1/\tilde{\nu}$, with $\tilde{\nu}=(X-e_3)/(e_1-e_3)$. Note that had we defined $\pi(X)$ as an integral from $X=\infty$ to $e_1$ instead of from $X = e_3$ to $e_2$ then it it would have been precisely equal to Jacobi's integral $\Pi_{\cal J}(\nu,K)$.

\subsection{Jacobi's elliptic integral of third kind in Weierstrass form}

Let us get a bit more specific than in (\ref{pi_W}) and define the elliptic integral of third kind in Weierstrass form by
\begin{equation}
\pi(X_1,X_2) = - \int_{\infty}^{X_1} \frac{Y_2}{X-X_2} \frac{dX}{2Y}  \label{pi(X1,X2)}
\end{equation}
In terms of Weierstrass elliptic functions, 
\begin{equation}
\pi(X_1,X_2) = \frac{1}{2} \ln \frac{\sigma(u_2+u_1)}{\sigma(u_2-u_1)} - u_1\zeta(u_2) \label{pi(X1,X2)-W}
\end{equation}
Alternatively, based on the representation formula for the Weierstrass $\sigma$-function in terms of Jacobi theta functions, one obtains the remarkably similar theta-function formula
\begin{equation}
\pi(X_1,X_2) = \frac{1}{2} \ln \frac{\vartheta_1(v_2+v_1,q)}{\vartheta_1(v_2-v_1,q)} - v_1\frac{\vartheta'_1(v_2,q)}{\vartheta_1(v_2,q)}
\end{equation}
with
\begin{equation}
v = \frac{\pi}{2} \frac{u}{\omega} \label{theta-vars}
\end{equation}
and the standard definition for $q$. Derivatives with respect to the first argument of the theta functions are denoted, conventionally, with primes. A more rigorous notation for the integral (\ref{pi(X1,X2)}) would be $\pi(u_1,u_2)$. Note for instance that, since $\sigma(u)$ and $\zeta(u)$ are odd functions, then so is $\pi(u_1,u_2)$ in either argument. Also, under lattice shifts, one has
\begin{eqnarray}
&& u_2 \longrightarrow u_2 + 2\omega_i \quad \Longrightarrow \quad \pi(X_1,X_2) \longrightarrow \pi(X_1,X_2) \\ [4pt]
&& u_1 \longrightarrow u_1 + 2\omega_i \quad \Longrightarrow \quad \pi(X_1,X_2) \longrightarrow \pi(X_1,X_2) + 2 \pi_i(X_2) \label{non-ell}
\end{eqnarray}
The notation $\pi(X_1,X_2)$ is completely obscure if not misleading with respect to properties such as (\ref{non-ell}). The careful reader is urged to use throughout $\pi(u_1,u_2)$ instead of $\pi(X_1,X_2)$, wherever the latter occurs in the text.

The following interchange of amplitude and parameter formula holds
\begin{equation}
\pi(X_1,X_2) - \pi(X_2,X_1) = 
\begin{array}{|cc|}
u_1		& u_2 \\
\zeta(u_1)	& \zeta(u_2) 
\end{array} 
+ \mathbb{Z} \cdot i\frac{\pi}{2}
\end{equation}

Jacobi's integrals satisfy an addition theorem. To see that, consider a set of $n+1$ points $(X_1,Y_1)$, $\cdots$, $(X_n,Y_n)$ and $(X,Y)$ on the Weierstrass curve together with the corresponding points $u_1$, $\cdots$, $u_n$ and $u$ on its Jacobian and note that from (\ref{pi(X1,X2)-W}) one has
\begin{equation}
\sum_{i=1}^{n} \pi(X_i,X) = \frac{1}{2} \ln \frac{F(-u,u_1,\cdots,u_n)}{F(\phantom{+}u,u_1,\cdots,u_n)} - \zeta(u) \sum_{i=1}^{n} u_i \label{lilo}
\end{equation}
with 
\begin{equation}
F(v,v_1, \cdots v_n) = \prod_{i=1}^n \frac{\sigma(v-v_i)}{\sigma(v)\sigma(v_i)}
\end{equation}
Now, provided that 
\begin{equation} \label{null-sum}
\sum_{i=1}^{n} u_i = 0
\end{equation}
one has 
\begin{equation}
F(u,u_1,\cdots,u_n) = \frac{1}{(n-1)!} \frac{\Delta_{(n)}(u,u_2,\cdots,u_n)}{\Delta_{(n-1)}(u_2,\cdots,u_n)} \label{Fu}
\end{equation}
with
\begin{equation}
\Delta_{(n)}(u_1,\cdots,u_n) = 
\begin{array}{|ccccc|}
1 & \wp(u_1) & \wp'(u_1) & \cdots & \wp^{(n-2)}(u_1) \\
 &  &  & &  \\
\vdots & \vdots & \vdots & & \vdots \\ 
 &  &  & &  \\
1 & \wp(u_n) & \wp'(u_n) & \cdots & \wp^{(n-2)}(u_n)
\end{array} \label{Delta_n}
\end{equation}
The proof of equation (\ref{Fu}) goes as follows \cite{MR1725862}: for fixed $u_1$, $\cdots$, $u_n$, both 
\begin{displaymath}
F(u,u_1,\cdots,u_n) \Delta_{(n-1)}(u_2,\cdots,u_n)
\qquad \mbox{and} \qquad
\Delta_{(n)}(u,u_2,\cdots,u_n)
\end{displaymath}
are meromorphic functions in $u$ with a pole of order $n$ at $u=0$ and simple poles at $u=u_1$, $\cdots$, $u_n$.  Their ratio is therefore a first order elliptic function, and hence a constant in $u$. To compute this constant one uses that $\sigma(u) = u + \cdots$ , $\wp(u) = 1/u^2 + \cdots$ and then compares the Laurent series of the two functions in a neighborhood of $u=0$.

In this case, equation (\ref{lilo}) becomes
\begin{equation}
\sum_{i=1}^{n} \pi(X_i,X) = \frac{1}{2} \ln \frac{\Delta_{(n)}(-u,u_2,\cdots,u_n)}{\Delta_{(n)}(\phantom{+}u,u_2,\cdots,u_n)} \label{lulu}
\end{equation}
As the derivatives of $\wp(u)$ are elliptic functions and so belong to the polynomial ring $\mathbb{C}[\wp,\wp']$, then, based on the fact that adding a multiple of a column to another leaves a determinant unchanged, one can show that (\ref{Delta_n}) is proportional to
\begin{equation}
\Xi_{(n)}(X_1,Y_1, \cdots  ,X_n,Y_n) = 
\begin{array}{|ccccccc|}
1 & X_1 & Y_1 & X_1^2 & Y_1X_1& X_1^3 & \quad \cdots \quad \\ 
\rule{0pt}{15pt} &&&&&& \\
\vdots &&&&& \vdots &  \\
\rule{0pt}{15pt} &&&&&& \\
1 & X_n & Y_n & X_n^2 & Y_nX_n & X_n^3 & \quad \cdots \quad
\end{array}
\end{equation}
Using this in equation (\ref{lulu}), one obtains the addition theorem for the Jacobi elliptic integrals of third kind
\begin{equation}
\sum_{i=1}^{n} \pi(X_i,X) = \frac{1}{2} \ln \frac{\Xi_{(n)}(X,-Y,X_2,Y_2, \cdots  ,X_n,Y_n)}{\Xi_{(n)}(X,\phantom{+}Y,X_2,Y_2, \cdots  ,X_n,Y_n)}
\end{equation} 

Note that (\ref{Fu}) implies also the addition formula for the Weierstrass $\wp$-function
\begin{equation}
\Delta_{(n)}(u_1,\cdots,u_n) = 0
\end{equation}
or, equivalently,
\begin{equation}
\Xi_{(n)}(X_1,Y_1, \cdots  ,X_n,Y_n) = 0
\end{equation}

Setting the upper limit in (\ref{pi(X1,X2)}) equal to $e_1$ yields the complete Jacobi integral of third kind, {\it i.e.},
\begin{equation}
\pi(e_1,X) = \pi(X) + \mathbb{Z} \cdot i\frac{\pi}{2}
\end{equation}
The ambiguity arises essentially because we define here $\pi(X)$ as an integral from $X = e_3$ to $e_2$ instead of from $X = \infty$ to $e_1$. In the process of deforming the second path into the first one picks a residue from the pole at $X$.

The complete integral has the following Weierstrass elliptic-function representation
\begin{equation}
\pi(X) = 
\begin{array}{|cc|}
u		& \omega \\
\zeta(u)	& \zeta(\omega) 
\end{array} \label{Jacobi_int}
\end{equation}
 and the theta-function representation
\begin{equation}
\pi(X) = -\frac{\pi}{2} \frac{\vartheta'_1(v,q)}{\vartheta_1(v,q)}   \label{piX-theta}
\end{equation}
with $v$ given in terms of $u$ as in (\ref{theta-vars}). This last relation is just the well-known  theta-function representation formula for the Weierstrass $\zeta$-function re-expressed in terms of $\pi(X)$.

An addition theorem holds also for the complete Jacobi integrals of third kind. Consider three points $u_1$, $u_2$, $u_3$ on the Jacobian, satisfying
\begin{equation}
u_1+u_2+u_3 = 0
\end{equation}
The corresponding points on on the Weierstrass curve are collinear
\begin{equation}
\begin{array}{|ccc|}
 \, 1  \,  	& X_1 	& \, Y_1 \,  \\ [5pt]
1 		& X_2	& Y_2  \\ [5pt]
1 		& X_3 	& Y_3
\end{array}
= 0
\end{equation}
The addition theorem for the Weierstrass $\zeta$-function together with the equation (\ref{Jacobi_int}) imply then that
\begin{equation}
\pi(X_1) + \pi(X_2) + \pi(X_3) = \omega \frac{Y_1-Y_2}{X_1-X_2} = \omega \frac{Y_2-Y_3}{X_2-X_3} = \omega \frac{Y_3-Y_1}{X_3-X_1}
\end{equation}
The ratios on the r.h.s. are equal (apart from the $\omega$ factor) to the slope of the line determined by the three colinear points under consideration. Interesting particular cases of this addition formula are obtained by setting for example $(X_3,Y_3) = (e_i,0)$ with $i = 1,2,3$, in which case $\pi(X_3)$ takes the values $0$, $ - i \pi/2$ and $+i \pi/2$, respectively. Such formulas have appeared in the nineteenth century mathematical literature in the guise of identities relating pairs of Legendre elliptic integrals of the third kind with characteristics $\nu$ satisfying special relations, see {\it e.g.} \cite{MR0124532}.

\subsection{Differentiation formulas} \label{ell-diff-formulas}

Legendre's complete elliptic integrals have the well-known differentiation formulas
\begin{eqnarray}
d K(k) &\hspace{-7pt} = \hspace{-7pt}&\frac{E(k)-k'^2 K(k)}{2k^2k'^2} dk^2 \label{dK} \\ [6pt]
d E(k) &\hspace{-7pt} = \hspace{-7pt}&\frac{E(k)-K(k)}{2k^2} dk^2 \label{dE} \\ [6pt]
d \Pi(\nu,k) &\hspace{-7pt} = \hspace{-7pt}&\frac{E(k)-k'^2 \Pi(\nu,k)}{2(k^2-\nu)k'^2} dk^2 \nonumber \\ 
\hspace{-7pt} &+& \hspace{-7pt} \frac{(\nu-k^2) K(k)-\nu E(k) - (\nu^2-k^2) \Pi(\nu,k)}{2\nu(\nu-1)(\nu-k^2)} d\nu \label{dPi}
\end{eqnarray}
Based on these, we derive the following differentiation formula for Jacobi's version of the complete elliptic integral of the third kind
\begin{equation}
d \Pi_{\cal J}(\nu,K) = -\frac{2\nu(\nu-1)dK+[(\nu-1)K+E] d\nu}{2\sqrt{\nu(\nu-1)(\nu-k^2)}} \label{dPi_J}
\end{equation}
This justifies the parametrization by $\nu$ and $K$.
Comparing with (\ref{dPi}), one can see that whereas the derivatives of Legendre's integral of third kind with respect to its parameters depend on the integral itself, the derivatives of Jacobi's do not. 

On the Weierstrass side one obtains, by means of the relations (\ref{K_JW}) through~(\ref{Pi_JW}),  as well as the equations (\ref{dK}), (\ref{dE}) and~(\ref{dPi_J}),  the following set of differentiation formulas 
\begin{eqnarray}
d\omega &\hspace{-7pt} = \hspace{-7pt}&-\frac{2g_2^2\omega-9 g_3\eta}{2\Delta} dg_2 + 3\frac{3g_3\omega-2g_2\eta}{2\Delta} dg_3 \label{d_omega} \\
d\eta &\hspace{-7pt} = \hspace{-7pt}&-g_2 \frac{3g_3\omega-2g_2\eta}{2\Delta} dg_2 + \frac{2g_2^2\omega-9 g_3\eta}{2\Delta} dg_3 \label{d_eta}
\end{eqnarray}
and
\begin{eqnarray}
d \pi(X) 	&\hspace{-7pt} = \hspace{-7pt}&\frac{\eta+X\omega}{2Y} dX \nonumber \\
		&\hspace{-7pt}  + \hspace{-7pt}& \frac{(g_2X+3g_3)(3g_3\omega-2g_2\eta) - X^2(2g_2^2\omega-9 g_3\eta)}{2Y\Delta}dg_2 \nonumber \\
		&\hspace{-7pt}  + \hspace{-7pt}& \frac{(3X^2-2g_2)(3g_3\omega-2g_2\eta) \,-\, X(2g_2^2\omega-9 g_3\eta)}{2Y\Delta} dg_3 \label{d_piX}
\end{eqnarray}
Note that the derivatives of $\pi(X)$ have again the remarkable feature that they do not depend on $\pi(X)$ itself, and are, moreover, elliptic functions. Occasionally we may use the shorthand notations
\begin{equation}
\eta^* = \frac{2g_2^2\omega-9 g_3\eta}{\Delta} \qquad\mbox{and}\qquad \omega^* = \frac{3g_3\omega-2g_2\eta}{\Delta} \label{omet1}
\end{equation} 
Incidentally, note that equations (\ref{d_omega}) and (\ref{d_eta}) constitute the first step in a recurrent series. For any integer $k \geq 0$ one has
\begin{eqnarray}
d\omega_k &\hspace{-7pt} = \hspace{-7pt}&\frac{1}{2} (-\eta_{k+1}dg_2 + [2\!+\!(-1)^k]\omega_{k+1}dg_3) \\
d\eta_k	&\hspace{-7pt} = \hspace{-7pt}&\frac{1}{2} (-g_2\omega_{k+1}dg_2 + [2\!-\!(-1)^k]\eta_{k+1}dg_3)
\end{eqnarray}
where
\begin{equation}
\eta_0 = \eta \qquad\quad \omega_0 = \omega
\end{equation}
and
\begin{equation}
\left(\!\! \begin{array}{c}
\omega_{k+1} \\ [4pt]
\eta_{k+1}
\end{array} \!\!\right)
=\frac{1}{\Delta} \!
\left(\!\! \begin{array}{cc}
[2\!-\!(-1)^k]\, 3g_3 	\!\!\!&\!\!\! 2g_2 \\ [4pt]
2g_2^2 			\!\!\!&\!\!\! [2\!+\!(-1)^k]\, 3g_3
\end{array} \!\!\right)
\!
\left(\!\! \begin{array}{cc}
6k+1 	\!&\!  0 \\ [4pt]
0 		\!&\! 6k-1
\end{array} \!\!\right)
\!
\left(\! \begin{array}{c}
\omega_k \\ [4pt]
\eta_k
\end{array} \!\right)
\end{equation}
This is proved by induction.

Observe that equations (\ref{d_omega}) and (\ref{d_eta}) can be reverted to express $dg_2$ and $dg_3$ in terms of $d\eta$ and $d\omega$, which amounts to regarding $g_2$ and $g_3$ as functions of $\eta$ and $\omega$ and not the other way around. From this perspective, $\pi(X)$ is then a function of $X$, $\eta$ and $\omega$. Its partial derivatives with respect to this alternative set of variables, as results from substituting the reverted equations (\ref{d_omega})-(\ref{d_eta}) into equation (\ref{d_piX}), are given by
\begin{equation}
d\pi(X) = \frac{\eta+X\omega}{2Y}dX - \frac{V\omega+X}{Y}d\eta - \frac{V\eta-X^2}{Y}d\omega \label{d_piX2}
\end{equation}
where
\begin{equation}
V = \frac{2g_2\eta - 3g_3\omega}{3\eta^2-g_2\omega^2} \label{vee}
\end{equation}

The Weierstrass roots $e_1$, $e_2$ and $e_3$ can be regarded as functions of the modular coefficients $g_2$ and $g_3$. As one can easily check by means of equations~(\ref{G2}) and (\ref{G3}), their partial derivatives with respect to these variables are given by
\begin{equation}
de_1 = \frac{e_1dg_2 + dg_3}{(e_1-e_2)(e_1-e_3)} \label{d_e123}
\end{equation}
and permutations thereof.

The Weierstrass roots can equally be regarded as functions of the period $\omega$ and the {\it nome} $q$. This dependence can be made in some sense explicit by means of the Jacobi theta-function representation formulas, which follow immediately from
\begin{eqnarray}
e_1-e_2 &\hspace{-7pt} = \hspace{-7pt}& \left(\frac{\pi}{2\omega}\right)^2\! \vartheta_4^4(0,q) \nonumber \\
e_1-e_3 &\hspace{-7pt} = \hspace{-7pt}& \left(\frac{\pi}{2\omega}\right)^2\! \vartheta_3^4(0,q)  \nonumber \\
e_2-e_3 &\hspace{-7pt} = \hspace{-7pt}& \left(\frac{\pi}{2\omega}\right)^2\! \vartheta_2^4(0,q) 
\label{W_roots_theta}
\end{eqnarray}
The partial derivatives of the roots with respect to $\omega$ are trivial to compute. The partial derivatives with respect to $q$ can be computed by using the fact that theta functions satisfy the 1-dimensional heat equation 
to exchange each derivative with respect to the elliptic modulus on a theta constant with a double derivative with respect the first argument of the corresponding theta function, evaluated at zero. By resorting to the equations~(\ref{W_roots_theta}) and the theta-function representation formulas (\ref{lola}), one can re-write the resulting theta-function expressions in terms of the roots. In the end, one obtains
\begin{equation}
d e_1 = -2\, e_1 \frac{d\omega}{\omega} + \frac{1}{3}\!\left(\frac{2}{\pi}\right)^2\!\omega\,[(e_1-e_2)(\eta+e_3\omega) + (e_1-e_3)(\eta+e_2\omega)] \frac{dq}{q} \label{dW_roots_dtau}
\end{equation}
and the permutations thereof.

\subsection{Theta-function representation formulas}

Differentiating equation (\ref{piX-theta}) with respect to $v$ and using the fact that $X = \wp(u) = -\zeta'(u)$, we get
\begin{equation}
\omega (\eta + X\omega) = -\left(\frac{\pi}{2}\right)^2 \frac{\partial^2}{\partial v^2} \ln \vartheta_1(v,q) \label{X-theta}
\end{equation}
Setting then in turns $X = e_1$, $e_2$, $e_3$ corresponding to which $v = \pi/2$, $-(1+\tau)\pi/2$, $\tau \pi/2$ respectively, and then making use of the well-known quasi-periodicity properties of the theta functions yields
\begin{equation}
\omega (\eta + e_{i\, }\omega) = - \left(\frac{\pi}{2}\right)^2 \frac{\vartheta_{i+1}^{''}(0,q)}{\vartheta_{i+1}(0,q)} \label{lola}
\end{equation}
with $i = 1,2,3$. Equations (\ref{X-theta}) and (\ref{lola}) can be interpreted as theta-function representation formulas for $X$ and the Weierstrass roots, respectively. For $\eta$ one has \cite{MR1054205}
\begin{equation}
\eta = -\frac{\pi^2}{12\omega} \frac{\vartheta_1'''(0,q)}{\,\vartheta_1'(0,q)} \label{eta-theta}
\end{equation}

The logarithmic derivatives of Jacobi's theta functions have the following Fourier series expansions (see {\it e.g.} \cite{Wolfram})
\begin{eqnarray}
\frac{\vartheta'_1(v,q)}{\vartheta_1(v,q)} &\hspace{-7pt} = \hspace{-7pt}& \phantom{+}\cot v + 4 \sum_{n=1}^{\infty} \frac{q^{2n}}{1-q^{2n}} \sin (2nv)  \label{theta1'/theta1} \nonumber \\
\frac{\vartheta'_2(v,q)}{\vartheta_2(v,q)} &\hspace{-7pt} = \hspace{-7pt}&  - \tan v + 4 \sum_{n=1}^{\infty} (-)^n \frac{q^{2n}}{1-q^{2n}} \sin (2nv) \nonumber \\
\frac{\vartheta'_3(v,q)}{\vartheta_3(v,q)} &\hspace{-7pt} = \hspace{-7pt}&  4 \sum_{n=1}^{\infty} (-)^n \frac{q^{n}}{1-q^{2n}} \sin (2nv) \nonumber \\
\frac{\vartheta'_4(v,q)}{\vartheta_4(v,q)} &\hspace{-7pt} = \hspace{-7pt}&  4 \sum_{n=1}^{\infty}  \frac{q^{n}}{1-q^{2n}} \sin (2nv) 
\label{Fourier}
\end{eqnarray}
By differentiating the last three relations with respect to $v$ and then setting $v = 0$ we obtain the $q$-series representations
\begin{eqnarray}
\frac{\vartheta''_2(0,q)}{\vartheta_2(0,q)} &\hspace{-7pt} = \hspace{-7pt}& -1 + 8 \sum_{n=1}^{\infty} (-)^n n \frac{q^{2n}}{1-q^{2n}} \nonumber \\
\frac{\vartheta''_3(0,q)}{\vartheta_3(0,q)} &\hspace{-7pt} = \hspace{-7pt}& 8 \sum_{n=1}^{\infty} (-)^n n \frac{q^{n}}{1-q^{2n}} \nonumber \\
\frac{\vartheta''_4(0,q)}{\vartheta_4(0,q)} &\hspace{-7pt} = \hspace{-7pt}& 8 \sum_{n=1}^{\infty}  n \frac{q^{n}}{1-q^{2n}}
\label{q-series}
\end{eqnarray}
The $q$-series expansion relevant for $\eta$ is \cite{MR1054205}
\begin{equation}
\frac{\vartheta_1'''(0,q)}{\vartheta_1'(0,q)} = -1 + 24 \sum_{n=1}^{\infty}n\frac{q^{2n}}{1-q^{2n}} \label{theta'''/theta'}
\end{equation}

\subsection{Modular transformations and $q'$-series expansions} \label{mod-transf}

The Jacobi theta-functions of {\it nome} $q=\exp(i\pi\tau)$ and those of complementary {\it nome} $q'=\exp(i\pi\tau')$ with 
\begin{equation}
\tau' = -\frac{1}{\tau} \label{tau-prime}
\end{equation}
are related as follows (see {\it e.g.} \cite{MR1054205})
\begin{eqnarray}
i \vartheta_1(v,q) &\hspace{-7pt} = \hspace{-7pt}& \sqrt{-i\tau'} e^{i\tau'v^2/\pi} \vartheta_1(\tau'v,q') \label{theta1-mod} \\ [4pt]
\vartheta_2(v,q) &\hspace{-7pt} = \hspace{-7pt}& \sqrt{-i\tau'} e^{i\tau'v^2/\pi} \vartheta_4(\tau'v,q') \\ [4pt]
\vartheta_3(v,q) &\hspace{-7pt} = \hspace{-7pt}& \sqrt{-i\tau'} e^{i\tau'v^2/\pi} \vartheta_3(\tau'v,q') \\ [4pt]
\vartheta_4(v,q) &\hspace{-7pt} = \hspace{-7pt}& \sqrt{-i\tau'} e^{i\tau'v^2/\pi} \vartheta_2(\tau'v,q') \label{theta4-mod}
\end{eqnarray}
This set of relations allows one to derive the modular transformation properties of any theta-function-dependent quantity. For instance, the Weierstrass coefficients $g_2$ and $g_3$ admit the following theta-function representation formulas
\begin{eqnarray}
g_2 &\hspace{-7pt} = \hspace{-7pt}& \hspace{3.5pt} \frac{1}{6}\, \left(\frac{\pi}{2\omega}\right)^4 \![\vartheta_2(0,q)^8\!+\vartheta_3(0,q)^8\!+\vartheta_4(0,q)^8] \label{g2-theta} \\
g_3 &\hspace{-7pt} = \hspace{-7pt}& \frac{1}{27} \left(\frac{\pi}{2\omega}\right)^6 \![\vartheta_2(0,q)^4\!+\vartheta_3(0,q)^4][\vartheta_3(0,q)^4\!+\vartheta_4(0,q)^4][\vartheta_4(0,q)^4\!-\vartheta_2(0,q)^4] \label{g3-theta}
\end{eqnarray}
Using the equations (\ref{theta1-mod})-(\ref{theta4-mod}) with $v$ set to 0 one can easily show that they remain invariant under the modular transformation (\ref{tau-prime}), {\it i.e.},
\begin{eqnarray}
g_2 &\hspace{-7pt} = \hspace{-7pt}& \hspace{3.5pt} \frac{1}{6}\, \left(\frac{\pi}{2\omega'}\right)^4 \![\vartheta_2(0,q')^8\!+\vartheta_3(0,q')^8\!+\vartheta_4(0,q')^8] \label{g2-theta-prime} \\
g_3 &\hspace{-7pt} = \hspace{-7pt}& \frac{1}{27} \left(\frac{\pi}{2\omega'}\right)^6 \![\vartheta_2(0,q')^4\!+\vartheta_3(0,q')^4][\vartheta_3(0,q')^4\!+\vartheta_4(0,q')^4][\vartheta_4(0,q')^4\!-\vartheta_2(0,q')^4] \label{g3-theta-prime}
\end{eqnarray}
On the other hand, quantities such as $\eta$ and $\pi(X)$ transform non-trivially under this modular transformation. Their theta-function representations are given by equations (\ref{eta-theta}) and (\ref{piX-theta}). The corresponding primed quantities have similar representations
\begin{eqnarray}
\eta' &\hspace{-7pt} = \hspace{-7pt}& -\frac{\pi^2}{12\omega'} \frac{\vartheta_1'''(0,q')}{\ \vartheta_1'(0,q')} \label{eta-prime-theta} \\ [4pt]
\pi'(X) &\hspace{-7pt} = \hspace{-7pt}& -\frac{\pi}{2} \frac{\vartheta'_1(\tau'v,q')}{\vartheta_1(\tau'v,q')}
\end{eqnarray}
and $q'$-expansion formulas given by the equations (\ref{theta'''/theta'}) and (\ref{theta1'/theta1}) with $q$ replaced by $q'$ and $v$ by $\tau' v$. Using the equations (\ref{theta1-mod})-(\ref{theta4-mod}) one can show then that the modular transformation rules for $\eta$ and $\pi(X)$ under (\ref{tau-prime}) are given precisely by the Legendre relations (\ref{Legendre-id}) respectively (\ref{Legendre-pi}). Their $q'$-series expansions follow immediately from these.

\vskip20pt
\noindent {\large \bf Acknowledgements} \\ [10pt]
The author wishes to thank Martin Ro\v{c}ek for his determining contributions to this project, for his constant encouragement as well as for reading through the manuscript.

\bibliographystyle{utphys}
\bibliography{article1}

\end{document}